\documentclass[11pt]{amsart}
\usepackage[francais]{babel}

\title{Particularit\'es des bases d'ordre 2 de $\N$\\ et  conjecture de Erd\"os-Tur\'an.} 

\author{Labib Haddad}
\address{120 rue de Charonne, 75011 Paris, France}
\email{labib.haddad@wanadoo.fr}

\usepackage{amssymb}
\usepackage{amsmath}
\usepackage{endnotes}

\newcommand{\su}{\subsection*}

\newcommand{\head}{\section*}
\newcommand{\noi}{\noindent}

\newcommand{\Ž}{\'e}

\newcommand{\ˆ}{\`a}
\newcommand{\}{\`u}

\newcommand{\N}{\mathbb N}

\newcommand{\Z}{\mathbb Z}

\newcommand{\cal}{\mathcal}

\newcommand{\leqs}{\leqslant}

\newcommand{\geqs}{\geqslant}

\newcommand{\nts}{\negthickspace}

\newcommand{\nms}{\negmedspace}

\newcommand{\ali} {\begin{aligned}}   
\newcommand{\ala} {\end{aligned}}

\newcommand {\et}{\ \text{et}\ }

\newcommand {\ou}{\ \text{ou}\ }

\newcommand {\si}{\ \text{si}\ }

\newcommand {\sinon}{\ \text{sinon}\ }

\newcommand {\pour}{\ \text{pour}\ }

\newcommand {\sia}{\ \text{if}\ }

\newcommand {\ora}{\ \text{or}\ }

\newcommand {\anda}{\ \text{and}\ }

\newcommand{\inc}{\subset}

\newcommand{\bc}{\begin{cases}}
\newcommand{\ec}{\end{cases}}

\begin{document}
\maketitle

\markboth{{\sc Labib Haddad}}{{\sc Bases de $\N$ et conjecture de Erd\"os-Tur\'an}}

\hfill{\it Tout peut na\"tre ici-bas d'une attente infinie.}

\hfill Paul Val\Žry

\

\su{R\Žsum\Ž}{\small On introduit un nouvel outil, {\sl l'alg\bre de Erd\"os-Tur\'an}, taill\Ž sur mesure pour l'\Žtude  des bases d'ordre 2 de $\N$. On verra comment cet outil est bien adapt\Ž \ˆ l'examen de la conjecture de Erd\"os-Tur\'an.  En des sens bien pr\Žcis, (topologique et probabiliste), on montre qu'il y a {\bf peu de place} et {\bf presque aucune chance} pour que la conjecture de Erd\šs-Tur\'an soit fausse. Sait-on jamais ?! 

\'Ecrit en fran\c cais, le texte de l'article est suivi de sa version en langue anglaise.} 

\su{Abstract}{\small A new tool is introduced,  specially tailored for the  study of order 2 bases of $\N$. We call it {the \sl Erd\"os-Tur\'an algebra} [or ring] since  it suits well the scrutiny of the Erd\"os-Tur\'an conjecture. In very precise meanings, (topological and measure theoretical) there is {\bf little room} and {\bf almost no chance} for the Erd\šs-Tur\'an conjecture to be wrong. Does one ever know\nms?\nms!

Written in French, the paper ends with its version in English.}

\

\head{Introduction}

\

\

\'Etant donn\Žs une  partie $X \inc \Z$ et un entier $z\in\Z$, d\Žsignons par $p(X,z)$ le nombre de tous les couples $(x,y)\in X \times X$ tels que $x+y=z$ et $x\leqs y$. Bien entendu,  $p(X,z)$ est \Žgal \ˆ un entier ou \ˆ l'infini, selon le cas.

\

\noi En 2005, parmi d'autres r\Žsultats bien plus g\Žn\Žraux, {\sc Nathanson} [9]  \Žtablit ceci : quelle que soit la famille d'entiers $(s_z)_{z\in \Z}$ o\ $s_z>0$, il existe une partie $X\inc \Z$ pour laquelle on a $p(X,z) = s_z$, pour chaque $z$.

\

\noi On se donne une partie $A\inc \N= \{0,1,2,\dots\}$ et, afin d'abr\Žger, on introduit le vocabulaire suivant. On appelle $A$-{\it repr\Žsentation} de $n$ tout couple $(x,y)\in A\times A$ tel que $x+y = n$ et on d\Žsigne par  $r_n$ le nombre de ces $A$-repr\Žsentations.  On dit que $A$ est une {\it base}  de $\N$ lorsque l'on a $r_n >0$, pour chaque $n\in \N$. S'il en est ainsi, on voit facilement que l'on doit avoir, n\Žcessairement, $0\in A \et 1\in A$, de sorte que l'on aura
$$r_0 = 1, r_1 = 2.$$
On peut voir aussi, sans grande difficult\Ž, que l'on a, n\Žcessairement,
$$r_2 = 1 \ou 3 \ , \ r_3 =r_2- 1 \ou r_2 + 1 \ , \ r_4 =  \frac{2r_3-r_2+1}{2}  \ou  \frac{2r_3- r_2+3}{2}.$$
On va montrer, plus g\Žn\Žralement, ceci.  Il existe une suite $d_n(b_2,\dots,b_n)$ de polyn\™mes   ayant la propri\Žt\Ž suivante : quelle que soit la base $A$,  et pour chaque $n\geqs1$, on a 
$$r_{n+1} = d_n(r_2,\dots,r_n) \si n+1\notin A,$$
$$ r_{n+1} = d_n(r_2,\dots,r_n) + 2 \si n+1\in A.$$
 Le degr\Ž de  ces polyn\™mes, en chacune des  variables, est  \Žgal \ˆ $1$ : ce  sont, ainsi, des fonctions multilin\Žaires des variables. Leurs coefficients sont dans l'anneau $\Z[1/2]$. Par exemple, et en particulier, on a
$$d_1 = 1, \ d_2(b _2) = b _2-1, \ 2d_3(b _2,b _3) = 2b _3 - b _2+1,$$
$$2d_4(b _2,b _3,b _4) = b _4  -\ 3b _3 - b _2 + b_3b_2 +  1,$$
$$2d_5(b _2,b _3,b _4,b _5) = 2b _5 - 3 b _4  + 5 b _3 + b _2+ b _4b _2-2  b _3b _2  - 1,$$
$$4d_6 = 4b _6 - 6 b _5 + 10 b _4 - 13 b _3+ 6 b _2 + 2b _5b _2  + 2 b _4 b _3 - 6 b _4 b _2 + 3 b _3 b _2   - 6.$$

\

\noi Cela contraste avec le r\Žsultat de {\sc Nathanson}, bien entendu, puisque la suite $(r_n)_{n\in \N}$, loin d'\tre quelconque, est assez {\bf contrainte} ! Insistons : cette suite ne se d\Žroule pas librement, {\it elle demeure toujours sous haute surveillance}.

\

\noi En 1941, {\sc Erd\šs} et {\sc Tur\'an} [2] se sont pos\Ž la question suivante : est-il vrai que la suite $(r_n)_{n\in \N}$ associ\Že \ˆ une base $A$ n'est jamais born\Že ? Ils ont \Žnonc\Ž la question sous la forme d'une conjecture qui revient \ˆ dire ceci.

\su{Conjecture \emph{ET}} {\sl Pour toute base $A$ de $\N$, la suite $(r_n)_{n\in \N}$ des nombres de $A$-repr\Žsentations n'est pas born\Že}.

\

\noi Cette conjecture peut s'\Žnoncer de diverses mani\res \Žquivalentes.  Parfois, l'\Žquivalence est  purement verbale, comme dans les deux exemples suivants. Il arrive cependant qu'elle soit plus profonde, comme on le verra, plus loin,  dans la suite.  Comme premier exemple, on peut \Žnoncer la conjecture sous la forme  suivante.

\su{Conjecture \emph{ET}}  \sl Pour toute base $A$ de $\N$, la suite des $r_n$ prend une infinit\Ž de valeurs distinctes\rm.

\su{Une dichotomie} Parmi les parties $A\inc \N$, une distinction entre deux classes s'impose, naturellement. Pour faire court, on dira que $A$ est  {\bf de classe inf\Žrieure} lorsque la suite associ\Že des nombres $r_n$ est born\Že. Sinon, on dira que $A$ est {\bf de classe sup\Žrieure}. 

\

\noi La conjecture de Erd\šs-Tur\'an s'\Žnonce ainsi de la mani\re  suivante.
\su{Conjecture \emph{ET}} {\sl Toute base  de $\N$ est de classe sup\Žrieure}.

\

\noi Pour \Žtudier la question, (voir, par exemple, [1], [3], [4], [6], [7], [8], [11]), on utilise un proc\Žd\Ž bien connu qui consiste \ˆ associer, \ˆ chaque partie $A\inc \N$,  deux s\Žries formelles, $f(t) \et g(t)$. La premi\re est la s\Žrie formelle repr\Žsentative de A :
$$f(t) = \sum_{a\in A} t^a.$$
Elle est de la forme
$$f(t) = \sum_{n\geqs 0} u_nt^n$$
o\ $u_n = 1 \ou 0$ suivant que $n$ appartient ou non \ˆ la partie $A$. La seconde est le carr\Ž de la premi\re, $g(t) = f(t)^2$. On v\Žrifie, sans grand d\Žtour, que l'on a
$$g(t) = \sum_{n\geqs 0} r_nt^n$$
o\ les $r_n$ sont pr\Žcis\Žment les nombres qui ont \Žt\Ž d\Žfinis ci-dessus. Si $A$ est une base, on doit avoir $u_0= u_1 = 1$, autrement dit, $0\in A$ et $1\in A$, de sorte que l'on a bien $r_0 = 1$ et $r_1 = 2$, comme on l'a dit. 

\

\noi On appellera $A$-{\it pr\Žsentation} de $n$ tout couple $(x,y)\in A\times A$ pour lequel on a $x+y= n \et x\leqs y$  et  on d\Žsignera par  $p_n$ le nombre de ces $A$-pr\Žsentations.  
\noi Afin d'\Žviter les p\Žriphrases, on dira que la suite $p = (p_0,p_1,p_2,\dots,p_n,\dots)$ est {\bf le profil} de la partie $A$. Bien qu'elles soient proches, et li\Žes, on ne doit  pas confondre $A$-pr\Žsentation et $A$-repr\Žsentation. De fait, on observe, sans d\Žtour, ceci.  On a $r_n = 2p_n - 1$, si $n/2$ appartient \ˆ $A$, et $r_n = 2p_n$, sinon. En particulier, pour $n$ impair, on a toujours $r_{n} = 2p_n$.  L'usage des $p_n$ se r\Žv\le parfois plus judicieux que celui des $r_n$, comme on le verra par la suite.

\

\noi Cette \Žtude  conduit, naturellement, \ˆ introduire l'alg\bre suivante que j'appellerai {\bf alg\bre de Erd\šs-Tur\'an} et que je d\Žsignerai par ET.  Voici comment on  la d\Žfinit.

\

\head{L'alg\bre de Erd\šs-Tur\'an}

\

Soit D $= \Z[1/2]$ l'anneau des nombres rationnels dyadiques, autrement dit, D$=  \{k/2^n : k\in \Z, n\in \N\}$. On se donne trois suites d'ind\Žtermin\Žes $(a_1, a_2, \dots)$, $(b_1, b_2, \dots)$, $(c_1, c_2,\dots)$. On d\Žsigne par L l'anneau des polyn\™mes en ces ind\Žtermin\Žes, \ˆ coefficients dans D, autrement dit 
$$\mathrm L = \mathrm D[a_1, a_2, \dots, b_1, b_2, \dots, c_1, c_2,\dots].$$
C'est l'alg\bre libre engendr\Že par les ind\Žtermin\Žes sur l'anneau D. On introduit  trois s\Žries formelles en $t$,  \ˆ coefficients dans L :
$$f(t) = 1 + a_1t + a_2t^2 + \dots = 1 + \sum_{n\geqs 1} a_nt^n,$$
$$g(t) = 1 + b_1t + b_2t^2+ \dots = 1 + \sum_{n\geqs 1} b_nt^n,$$
$$h(t) = 1 + c_1t + c_2t^2+ \dots = 1 + \sum_{n\geqs 1} c_nt^n.$$
On impose les \Žgalit\Žs suivantes lesquelles \Žtablissent des liens entre les ind\Žtermin\Žes :
\[g(t) = f(t)^2 \et 2h(t) = g(t) + f(t^2).\tag{0}\]
On obtient un syst\me infini de relations entre ces ind\Žtermin\Žes.
On a ainsi, plus pr\Žcis\Žment,  deux types de relations. En posant $a_0=1$ et en calculant le carr\Ž de la s\Žrie $f(t)$, on obtient un premier type de relations que voici :
\[b_n = a_0a_n + a_1a_{n-1} + \dots + a_ka_{n-k} + \dots +
a_na_0.\tag{$*$}\]
Les relations du second type proviennent de la seconde identit\Ž, celle qui lie $f, g \et h$. Elles s'\Žcrivent :
\[2c_{2n+1} = b_{2n+1} \et 2c_{2n} = b_{2n} + a_n.\tag{$**$}\]
L'alg\bre ET est alors, {\it par d\Žfinition}, l'alg\bre L soumise \ˆ l'ensemble des  relations  ($*$) et ($**$).

\

\noi Dans les formules ($**$), on peut regrouper les deux cas, pair et  impair, en un seul. Par convention, on pose  $a_x = 0$, pour tout $x\notin \N$, et l'on a :
\[2c_n = b_n + a_{n/2}.\tag{$**$}\]
Les formules ($*$) sont simples et expriment les $b_n$ en fonction des $a_n$. Il existe des formules semblables, moins simples, pour exprimer les $a_n$ en fonction des $b_n$.  {\it Plus cach\Žes}, elles m\lent {\bf nombres de Catalan} et {\bf polyn\™mes de Bell} ordinaires. Pour les obtenir, on utilise les ressources fournies par les racines carr\Žes de s\Žries formelles. On fera cela, pas \ˆ pas, en donnant les d\Žtails n\Žcessaires.

\

\head{Racines carr\Žes de s\Žries formelles}

\su{A1 Nombres de Catalan} Les nombres de Catalan,
$C_1,C_2,\dots,C_k,\dots,$ sont d\Žfinis comme suit :
$$C_1 = 1\ \et \ C_k = C_1C_{k-1} + C_2C_{k-2} +
\dots C_{k-1}C_1, \pour k \geq 2.$$ 
Ce sont donc des {\bf entiers positifs}.  En voici les premi\res valeurs :
$$C_1 = 1 \ , \ C_2 = 1 \ , \ C_3 = 2 \ , \ C_4 = 5 \ , \
C_5 = 14.$$
En vertu de la relation de r\Žcurrence, leur s\Žrie g\Žn\Žratrice, la s\Žrie formelle $C(t) = \sum_{k \geqs 1} C_kt^k$,
satisfait   l'identit\Ž $C = t + C^2$.  On a, ainsi
$$C(t) = \frac{1 - \sqrt{1 - 4t}}{2} \ , \ \text{de sorte que}
\
\ \sqrt{1 - 4t} = 1
 - 2\sum_{k \geqs 1} C_k t^k.$$ 
\`A l'aide de la formule du bin\™me, on a
$$ \sqrt{1 - 4t} = 1 + \sum_{k
\geqs 1}
\binom {1/2}{k}(-4)^kt^k.$$
D'o\
$$C_k = \frac{(-1)^{k-1}}{2} 4^k \binom {1/2}{k} =
\frac{1}{k}\binom {2k - 2}{k -1}.$$ 
Par la formule de Wallis, il vient 
$$\lim_{n \to \infty} 4 \frac{C_k}{4^k} k^{3/2} =
\frac{1}{\sqrt\pi}.$$

\

\su{A2  Polyn\™mes de Bell} Les polyn\™mes de Bell
(ordinaires)   sont les polyn\™mes $P_{n,k}$ d\Žfinis par les relations  suivantes :
\[(b_1t + \dots + b_nt^n + \cdots)^k =  P_{k,k}t^k +
P_{k+1,k}t^{k+1} + \dots + P_{n,k}t^n +
\cdots,\tag {1}\]
pour des entiers $n \geqs k
\geqs 0$. Dans les autres cas, on convient de poser
$P_{n,k} =0$. 

\

\noi  Chacun des $P_{n,k}$ est un polyn\™me en les variables
$b_1,b_2,\dots,b_{n-k+1}$ dont tous les coefficients sont des entiers positifs, comme il appara\"t clairement au travers des relations qui les d\Žfinissent. Par souci de simplicit\Ž, on pourra abr\Žger l'\Žcriture $P_{n,k}(b_1,b_2,\dots,b_{n-k+1})$
  en
$P_{n,k}[b]$ o\ l'on convient avoir pos\Ž $b = b_1t + \dots + b_nt^n + \cdots$.

\

\noi Il faut se souvenir que, par convention, on a pos\Ž $P_{n,k} = 0$ lorsque les in\Žgalit\Žs $n \geqs k \geqs 0$ ne tiennent pas ! Cela apporte de grandes simplifications dans le maniement
 des indices de sommation.

\su{A3 Racines carr\Žes de s\Žries formelles} On reprend la s\Žrie formelle \ˆ coefficients dans l'alg\bre ET :
$$f(t) = 1 + a_1t + \dots + a_nt^n +
\cdots,$$  
et  son carr\Ž :
$$g(t) = f(t)^2 = 1 + b_1t + \dots + b_nt^n + \cdots.$$
On en a d\Žduit  les formules simples suivantes :
\[b_n = a_0a_n + a_1a_{n-1} + \dots + a_ka_{n-k} + \dots +
a_na_0.\]
Voici les formules qui redonnent les
coefficients
$a_n$ comme fonctions des  $b_k$ : 
\[a_n = 2 \sum_{1 \leqs k \leqs n}
\frac{(-1)^{k-1}}{4^k} C_k
P_{n,k}(b_1,b_2,\dots,b_{n-k+1}).\tag{2}\]

\su{D\Žmonstration} On pose, simplement, $b = b(t) = b_1t +
\dots + b_nt^n + \cdots,$ et on utilise la formule du bin\™me pour avoir
$$f(t) = (1+b)^{1/2} = 1 + \sum_{k \geqs 1} \binom{1/2}{k}
b^k.$$ 
Ensuite, \sl \ˆ  la Faa di Bruno\rm, on substitue \ˆ
$b^k$ son d\Žveloppement  donn\Ž par la formule (1), et l'on obtient
$$1 + \sum_{n \geqs 1} a_nt^n = 1 + 2 \sum_{n \geqs 1} \sum_{k
= 1}^n
\frac{(-1)^{k-1}}{4^k}C_k P_{n,k}[b]t^n.$$
Les formules (2) s'obtiennent en identifiant les coefficients de
$t^n$, des deux c\™t\Žs de l'identit\Ž.\qed

\

\noi Pour davantage de details, voir, par exemple,  {\sc Louis COMTET}, {\it Analyse combinatoire} I, P.U.F., Paris, 1970.

\su{A4 Les multi-indices} On utilisera les {\bf multi-indices} dont l'usage est particuli\rement indiqu\Ž pour traiter des polyn\™mes de  Bell. Un multi-indice, plus pr\Žcis\Žment, un
{\bf $k$-indice} est une 
$k$-sequence  (i.e., une suite finie de $k$ termes) $\mu =
(\mu_1,\dots,\mu_k)$ d'entiers $\mu_j > 0$.
\`A un tel multi-indice est associ\Ž son
\bf poids\rm, l'entier
$|\mu| = \mu_1 + \mu_2 +\dots + \mu_k$. Une s\Žrie
$b = b_1t + b_2t^2 + \dots + b_nt^n + \cdots$ \Žtant donn\Že, on associe \Žgalement au  $k$-indice
$\mu = (\mu_1,\dots,\mu_k)$ le {\bf mon\™me}
$b_\mu = b_{\mu_1}b_{\mu_2}\dots b_{\mu_k}$. On d\Žsigne par $M(k,n)$ l'ensemble des $k$-indices dont le poids est \Žgal \ˆ
$n$. Alors, comme on le voit sans d\Žtour \ˆ partir des relations qui les d\Žfinissent, les polyn\™mes de Bell  peuvent s'\Žcrire comme ceci :
$$P_{n,k}[b] = \sum_{ \mu \in M(k,n)}
b_\mu.$$
On a, par exemple, 
$$P_{n,1}[b] = b_n \ , \  P_{n,2}[b] = b_1b_{n-1} + b_2
b_{n-2} + \dots + b_{n-1} b_1 \ , \   P_{n,n}[b] = b_1^n.$$
D'apr\s les relations ($*$), lorsque les $a_n$ sont des entiers, il en est de m\me des 
$b_n$  donc, \Žgalement, des
$b_\mu$.

\su{A5 Coefficients de la s\Žrie racine carr\Že} La formule (2) s'\Žcrit alors
$$2a_n = 2a_n[b] = 4\sum_{\mu \in M(k,n)}
\frac{(-1)^{k-1}}{4^k} C_k  b_\mu.$$  
En particulier, il vient

\

$2a_1  = b_1$,

$2a_2  = b_2 - \frac{1}{4}b_1^2$,

$2a_3  = b_3 - \frac{1}{2} b_1b_2 +
\frac{1}{8}b_1^3$,

$2a_4 = b_4 - \frac{1}{2}b_1b_3 - \frac{1}{4}b_2^2 + \frac{3}{8}b_1^2b_2 - \frac{5}{64}b_1^4$,

$2a_5 = b_5 -\frac{1}{2}b_1b_4 - \frac{1}{2}b_2b_3 + \frac{3}{8}b_1^2b_3 + \frac{3}{8}b_1b_2^2 - \frac{5}{16}b_1^3b_2 + \frac{7}{128}b_1^5$,

$2a_6  = b_6 - \frac{1}{2}b_1b_5 - \frac{1}{2}b_2b_4 + \frac{3}{8}b_1^2b_4 - \frac{5}{16}b_1^3b_3- \frac{1}{4}b_3^2 -\frac{15}{32}b_1^2b_2^2  +  \frac{3}{4}b_1b_2b_3$

$\qquad \ \ + \frac{1}{8}b_2^3 + \frac{35}{128}b_1^4b_2- \frac{21}{512}b_1^6$,

$\dots\dots\dots\dots\dots\dots\dots$

\[2a_n  = b_n - \frac{1}{4}(b_1b_{n-1} + b_2
b_{n-2} + \dots + b_{n-1} b_1)
+ 4\sum_{\underset{ k > 2}{\mu\in M(k,n)}}
\frac{(-1)^{k-1}}{4^k} C_k  b_\mu.\tag{3}\]

\

\noi Il en d\Žcoule la relation remarquable suivante sur laquelle on reviendra :
\[b_n = 2a_n + \frac{1}{4}(b_1b_{n-1} + b_2
b_{n-2} + \dots + b_{n-1} b_1)
+ 4\sum_{\underset{ k > 2}{\mu\in M(k,n)}} \frac{(-1)^k}{4^k}C_k  b_\mu.\tag{4}\]

\

\noi On observe que le  terme
$\sum_{\underset{ k > 2}{\mu\in M(k,n)}}
\frac{(-1)^{k-1}}{4^k} C_k  b_\mu$
est un polyn\™me qui d\Žpend de $b_1,\dots,b_{n-2}$ mais pas de $b_{n-1}$ : en effet, pour  $k>2$ et pour chacun des $k$-indices $\mu= (\mu_1,\dots,\mu_k)$ de poids $n$,  on a n\Žcessairement $\mu_j \leqs n-2$. On \Žcrira ainsi
$$\sum_{\underset{ k > 2}{\mu\in M(k,n)}}
\frac{(-1)^{k-1}}{4^k}C_k  b_\mu = Q_n(b_1,\dots,b_{n-2}).$$  

\su{En passant} On remarquera  ceci, au passage. La conjecture {\bf\emph{ET}} de Erd\šs-Tur\`an revient \ˆ dire que, pour toute suite {\it born\Že} d'entiers $b_k\geqs 1$, il existe un indice $n$ pour lequel le $a_n$ de la formule (3)  n'est ni $0$ ni $1$.

\

\su{A6 Un cas particulier} \sl Le cas o\ chacun des $a_n$ est \Žgal \ˆ $0$ ou \ˆ $1$\rm.

\noi Dans ce cas, on aura $a_n^2 = a_n$ et cette relation implique la suivante :
$$\left(\frac{1}{2}b_n - \frac{1}{8}(b_1b_{n-1} + b_2
b_{n-2} + \dots + b_{n-1} b_1)
+ 2Q_n(b_1,\dots,b_{n-2})\right)^2 =$$
$$\frac{1}{2}b_n - \frac{1}{8}(b_1b_{n-1} + b_2
b_{n-2} + \dots + b_{n-1} b_1)
+ 2Q_n(b_1,\dots,b_{n-2}).$$ 
On en d\Žduit  une expression de $b_n^2$ comme polyn\™me en  $b_1, b_2, \dots, b_n$, dont le degr\Ž  en $b_n$ est \Žgal \ˆ $1$.
En partant de
$$b_{n+1} = 2a_{n+1} + \frac{1}{4}(b_1b_{n} + b_2
b_{n-1} + \dots + b_{n} b_1) + 4 Q_{n+1}(b_1,\dots,b_{n-2},b_{n-1}),$$
on y  remplace, {\it successivement},  $b_{n-1}^2, b_{n-2}^2,\dots$ par leurs valeurs respectives. On obtient une expression de $b_{n+1}$ sous la forme
$$b_{n+1} = 2a_{n+1} + d_n$$
o\ $d_n$ est un polyn\™me en $b_1, b_2, \dots, b_n,$ dans lequel chacune de ces variables figure avec un degr\Ž au plus \Žgal \ˆ $1$.

\

\noi Cela sugg\re une {\it sp\Žcialisation} de l'alg\bre ET que voici.

\

\head{L'alg\bre de Erd\šs-Tur\'an sp\Žciale}

\

J'appelle {\bf alg\bre de Erd\šs-Tur\'an sp\Žciale}, et je la d\Žsigne par ETS, l'alg\bre que l'on obtient \ˆ partir de ET  en adjoignant la suite des relations d'idempotence suivantes :
\[a_n^2 = a_n, \pour n = 1, 2, \cdots.\tag{$***$}\]
Autrement dit, ETS s'obtient \ˆ partir de l'alg\bre libre L en lui adjoignant l'ensemble des relations ($*$), ($**$) et ($***$).

\

\noi Bien entendu, de ces relations fondamentales, de nombreuses autres d\Žrivent.

\

\noi Par exemple, en rapprochant ($**$) de ($*$), on obtient les identit\Žs suivantes :
$$c_{2n} = a_{2n} +  a_{2n-1} + a_2a_{2n-2} + \dots + a_{n-1}a_{n+1} + a_n,$$
$$c_{2n+1} =  a_{2n+1} +  a_{2n} + a_2a_{2n-1} + a_3a_{2n-2}+\dots + a_na_{n+1},$$
Avec la convention $a_0=1$, on peut les regouper en une seule que voici :
\[c_n = \sum_{0\leqs k\leqs n/2} a_ka_{n-k}.\tag{5}\]

\

\noi Notamment, aussi, du rapprochement de (3) et ($***$) d\Žcoule  une expression de la forme
\[b_{n+1} = 2a_{n+1} + d_n(b_1,b_2,\dots,b_n)\tag{6}\]
o\ $d_n$ est un polyn\™me en $b_1, b_2, \dots, b_n,$ dans lequel chaque ind\Žtermin\Že figure avec un degr\Ž au plus \Žgal \ˆ $1$.
Comme dans le cas particulier trait\Ž ci-dessus, on a un algorithme (simple) qui permet de calculer chacun des polyn\™mes $d_n$. 

\

\su{Un  algorithme g\Žn\Žral de r\Žduction} C'est un proc\Žd\Ž d'{\bf \Žlimination}. On se donne une suite de variables $x_1,x_2,\dots$, li\Žes par des relations de la forme
$$x_n^2 = t_n(x_1,\dots,x_n)$$
o\ $t_n(x_1,\dots,x_n)$ est un polyn\™me en $x_1,\dots,x_n,$ dont le degr\Ž en $x_n$ est au plus \Žgal \ˆ 1. Tout polyn\™me $P(x_1,\dots,x_n)$ se r\Žduit alors \ˆ un polyn\™me en $x_1,\dots,x_n,$ o\ chacune des variables figure avec un degr\Ž au plus \Žgal \ˆ 1, par \Žlimination des puissances $x_k^j$, pour $j\geqs 2$.

\

\noi En effet, pour chaque entier $j\geqs 2$, on a une expression de $x_n^j$ comme polyn\™me en $x_1,\dots,x_n,$ dont le degr\Ž en $x_n$ est au plus \Žgal \ˆ 1. Cela se fait par r\Žcurrence, en commen\c cant par $x_n^3 = t_n(x_1,\dots,x_n)x_n$. 

\

\noi On remplace dans le polyn\™me $P(x_1,\dots,x_n)$ chacune des puissances $x_n^j$ par son expression et l'on obtient un polyn\™me $P_n(x_1,\dots,x_n)$ o\ $x_n$ figure avec le degr\Ž 1, au plus. En r\Žit\Žrant le proc\Žd\Ž pour $x_{n-1}, x_{n-2},\dots, x_1$, on obtient les polyn\™mes $P_{n, n-1}, P_{n, n-1, n-2}, \dots, P_{n, n-1,\dots,1}$,  et le r\Žsultat annonc\Ž.

\su{Une variante} Dans l'algorithme pr\Žc\Ždent, au lieu de remplacer dans le polyn\™me $P$ chacune des puissances $x_n^j$ par son expression, on peut utiliser, simplement,  la division euclidienne. Plus pr\Žcis\Žment,  on effectue la division du polyn\™me $P(x_1,\dots,x_n)$ par le polyn\™me $T(x_n) = x_n^2 - t_n(x_1,\dots,x_n)$ suivant la variable $x_n$ :
$$P(x_n) = Q(x_n)T(x_n) + R(x_n).$$
Le reste $R(x_n)$ ainsi obtenu est, pr\Žcis\Žment, le polyn\™me $P_n$. On fait cela \ˆ chacune des \Žtapes $n-1,n-2,\dots,1$.

\su{La tour d'alg\bres} On a ainsi un \Žtagement  d'alg\bres de polyn\™mes :
$$\mathrm L = \mathrm D[a_1, a_2, \dots, b_1, b_2, \dots, c_1, c_2,\dots]$$
$$\mathrm{ET   = L}\  \  \mathrm{modulo} \ \ (*) \et (**)$$
$$\mathrm{ETS = ET} \  \ \mathrm{modulo} \ \ (***)$$

On va lui ajouter un nouvel \Žtage, {\it pour les bases}, par sp\Žcialisation :
$$\mathrm{ETB = ETS}\  \ \mathrm{modulo}\ \  a_1 = 1.$$

\head{Le cas des bases\\
l'alg\bre ETB}

\

Dans le cas des bases, on a toujours $r_1 = 2$. On sp\Žcialise  alors l'alg\bre ETS en une alg\bre ETB, en posant $a_1 = 1$, $b_1=2$ : c'est un outil bien  {\bf adapt\Ž} \ˆ l'\Žtude des nombres $r_n$  et $p_n$ de repr\Žsentations et de pr\Žsentations associ\Žs aux bases.  

\

\centerline{\{Partout dans la suite, on adopte la convention $a_0 = a_1 = c_0=c_1 =1$.\}}

\

\noi Dans cette alg\bre, on a de nombreuses relations d\Žriv\Žes int\Žressantes dont voici les premi\res, des cas particuliers de (3) :

\

$2a_2  = b_2 - 1$,

$2a_3  = b_3 -  b_2 + 1$,

$2a_4 = b_4 - b_3 - \frac{1}{4}b_2^2 + \frac{3}{2}b_2 - \frac{5}{4}$,

$2a_5 = b_5 -b_4 - \frac{1}{2}b_2b_3 + \frac{3}{2}b_3 + \frac{3}{4}b_2^2 - \frac{5}{2}b_2 + \frac{7}{4}$,

$2a_6 = b_6 - b_5 - \frac{1}{2}b_2b_4 + \frac{3}{2}b_4 - \frac{5}{2}b_3+ \frac{35}{8}b_2 - \frac{1}{4}b_3^2 -\frac{15}{8}b_2^2  +  \frac{3}{2}b_2b_3+ \frac{1}{8}b_2^3 - \frac{21}{8}$.

\

\noi On \Žcrit $a_2^2 = a_2$. Il vient $b_2^2 = 4b_2 -3$. En rempla\c cant $b_2^2$ par cette valeur dans les expressions de $2a_4$ et $2a_5$, on obtient

\

 $2a_4 = b_4 - b_3 - \frac{1}{4}(4b_2-3) + \frac{3}{2}b_2 - \frac{5}{4}= b_4 - b_3+\frac{1}{2}b_2 -\frac{1}{2}$.
 
 $2a_5 = b_5 -b_4 - \frac{1}{2}b_2b_3 + \frac{3}{2}b_3 + \frac{3}{4}(4b_2 - 3)- \frac{5}{2}b_2 + \frac{7}{4}= b_5 -b_4 -\frac{1}{2}b_2b_3+$
 
  \qquad \ \ $\frac{3}{2}b_3 +  \frac{1}{2}b_2- \frac{1}{2} $.

\

\noi On \Žcrit $a_3^2 = a_3$. Il vient $b_3^2 = 2b_2b_3 -4b_2 + 4$, $b_2^3 = 13b_2 - 12$. On remplace $b_2^2$, $b_2^3$ et $b_3^2$ par leurs valeurs respectives dans l'expression de $2a_6$ et on obtient

\

$2a_6 = b_6 - b_5 + \frac{3}{2}b_4 - \frac{1}{2}b_2b_4  - \frac{5}{2}b_3+ \frac{35}{8}b_2 - \frac{1}{4}(2b_2b_3 -4b_2 + 4)-\frac{15}{8}(4b_2-3)  +  \frac{3}{2}b_2b_3+ \frac{1}{8}(13b_2 -12) - \frac{21}{8}= b_6 - b_5 + \frac{3}{2}b_4 - \frac{5}{2}b_3 + 24b_2 + b_2b_3 -\frac{1}{2}b_2b_4-\frac{1}{2}= $.

$b_6 -b_5 + \frac{3}{2}b_4 - \frac{5}{2}b_3 - \frac{1}{2}b_2 + b_2b_3 - \frac{1}{2}b_2b_4 + \frac{1}{2}$.

\

\noi On peut  faire les calculs {\it \ˆ la main}; cela va plus vite sur ordinateur. On retrouve les formules annonc\Žes dans l'introduction et, en particulier,
$$2d_5(b _2,b _3,b _4,b _5) = 2b _5 - 3 b _4  + 5 b _3 + b _2+ b _4b _2 - 2 b _3b _2  - 1,$$
$$4d_6 = 4b _6 - 6 b _5+ 10 b _4 - 13 b _3+ 6 b _2 + 2b _5b _2  + 2 b _4b _3  - 6 b _4 b _2 + 3 b _3 b _2   - 6.$$

\

\noi On fait un pas de plus en exprimant ces formules \ˆ l'aide des $c_n$ au lieu des $b_n$, en utilisant les relations :
\[2c_{2n+1} = b_{2n+1} \et 2c_{2n} = b_{2n} + a_n.\tag{$**$}\]
Ainsi, 
\

$2a_2  = b_2 - 1$ 

$2a_3  = b_3 -  b_2 + 1$ 

$2a_4 =  b_4 - b_3+\frac{1}{2}b_2 -\frac{1}{2}$
  
$2a_5 = b_5 -b_4 -\frac{1}{2}b_2b_3+\frac{3}{2}b_3 +  \frac{1}{2}b_2- \frac{1}{2}$ 

\noi deviennent, respectivement,

$a_2 = c_2 -1$,

$a_3 = c_3 - c_2 + 1$,

$a_4 = c_4 - c_3$,
 
$a_5  = c_5 - c_4 + 2c_3 + c_2 - c_2c_3 - 1$.

\

{\bf On voit par l\ˆ que ces expressions des $\bf a_n$ en fonction des $\bf c_k$ sont plus simples que les expressions pr\Žc\Ždentes en fonction des $\bf b_k$.} En effet, ce sont des polyn\™mes en $c_k$ dont tous les coefficients sont entiers.

\

\noi On peut ainsi se poser la question :  $a_n$ s'exprime-t-il toujours comme un polyn\™me de l'anneau $\Z[c_2,\dots,c_n]$ ? La r\Žponse est oui ! On va le montrer, dans un instant, apr\s avoir \Žtabli le r\Žsultat plus g\Žn\Žral suivant.

\

\head{Formes r\Žciproques}

\

Afin de faire court, dans toute la suite du texte, on utilisera le n\Žologisme  suivant : on qualifiera de  {\bf compliforme} tout polyn\™me  \ˆ coefficients entiers, de degr\Ž $1$ (au plus) en chacune des variables. [Ce n\Žologisme n'a qu'un statut \Žph\Žm\re, jusqu'\ˆ la fin du texte, aucune pr\Žtention ˆ la p\Žr\Žnit\Ž.] 

\

\noi De m\me, pour abr\Žger les \Žcritures, et tout au long du texte, on d\Žsignera par M l'alg\bre des polyn\™mes en $x_2,x_3,\dots$, \ˆ coefficients entiers, et par C l'ensemble des polyn\™mes compliformes en $x_2,x_3,\dots$., de sorte que l'on a ainsi 
C $\inc$ M $= \Z[x_2,x_3,\dots]$.

\

\su{Proposition} \sl Dans un anneau commutatif $R$ quelconque, on se donne deux suites, $y_2,y_3,\dots,$ et $z_2,z_3,\dots,$ o\ les $y_n$ sont des idempotents,  et telles que l'on ait :
$$z_{m+1} = y_{m+1} + s_m(y_2,\dots,y_m)$$
o\ $s_m\in$ M\rm.

\noi (i) \sl On aura
$$y_{n+1} = z_{n+1}  + t_n(z_2,\dots,z_n)$$
o\  $t_n\in$ C, {\rm i.e.,} les $t_n$ sont des polyn\™mes compliformes\rm.

\

\noi (ii) \sl De plus, pour chaque polyn\™me  $u(x_2,\dots,x_n)\in$ M, il existe un polyn\™me $v(x_2,\dots,x_n)\in$ C tel que l'on ait 
$$u(y_2,\dots,y_n)= v(z_2,\dots,z_n).$$\rm

\su{D\Žmonstration}  On proc\de par r\Žcurrence. On suppose  que, pour un $m\geqs 2$ donn\Ž et pour chaque indice $k$, $2\leqs k \leqs m$,  l'\Žl\Žment $y_k$ poss\de une expression sous  la forme  $y_k = z_k + t_{k-1}(z_2,\dots,z_{k-1})$ o\ $t_{k-1}$ est un polyn\™me compliforme. Cela est vrai pour $m = 2$ car $y_2 = z_2 - s_1$ et $s_1$ est un entier.  On \Žcrit $y_k^2 = y_k$ et on obtient une expression de $z_k^2$   comme polyn\™me \ˆ coefficients entiers en $z_2,\dots,z_k$, de degr\Ž $1$ en $z_k$. \`A l'aide de l'algorithme ci-dessus, on peut {\bf r\Žduire}  tout polyn\™me \ˆ coefficients entiers, en $z_2,\dots,z_k$, \ˆ un polyn\™me compliforme. On proc\de \ˆ la r\Žduction du polyn\™me
$$ s_m(\dots, z_k +  t_{k-1}(z_2,\dots,z_{k-1}), \dots)= s_m(y_2,\dots,y_m)= z_{m+1} - y_{m+1}.$$
L'expression obtenue ach\ve la r\Žcurrence, et la d\Žmonstration de la premi\re partie (i). Pour le (ii),  on \Žcrit
$$u(y_2,\dots,y_n) = u(z_2 + t_1, \dots, z_n + t_{n-1}(z_2,\dots,z_n))$$
et on r\Žduit ce dernier polyn\™me  \ˆ un polyn\™me  compliforme.
\qed

\

\su{Corollaire}{\sl Pour chaque entier $n\geqs 1$, il existe un polyn\™me $e_n(x_2,\dots,x_n)$ compliforme tel que, dans l'alg\bre ETB, on ait :
\[c_{n+1} = a_{n+1} + e_n(c_2,\dots,c_n)\tag{7}\]}
De plus, il  existe  un algorithme simple pour calculer ces polyn\™mes.

\

\noi En effet, il suffit d'utiliser la proposition pr\Žc\Ždente, en \Žcrivant l'identit\Ž
\[c_n = \sum_{0\leqs k\leqs n/2} a_ka_{n-k},\tag{5}\]
sous la forme
$$c_{m+1} = a_{m+1} + s_m(a_2,\dots,a_m).$$

\

\noi La proposition ne s'applique  pas aux   $b_n$ car
$b_n = 2a_n + \sum_{1\leqs k < n} a_ka_{n-k}$ n'est pas de la forme requise.

\

\su{Un premier cas tr\s particulier} En poursuivant les calculs pr\Žc\Ždents de $a_2,a_3,a_4,a_5$, on arrive \ˆ la formule suivante :
$$a_6 = c_6 - c_5 + 2c_4 - 4c_3 + 2c_3c_2 - c_4c_2.$$
Elle rec\le une petite surprise : si l'on y fait $c_2 = c_3 = c_4 = c_5 = c_6 = 1$, on arrive \ˆ $a_6 = -1$ qui est incompatible avec une base ! On peut \Žnoncer ce r\Žsultat comme suit : soient $p_2, p_3, p_4, p_5,p_6$ les nombres  de $A$-pr\Žsentations associ\Žs \ˆ une base donn\Že $A$ quelconque; alors  au moins l'un d'eux est $\geqs 2$ !

\noi Cela implique, bien entendu, que l'un au moins de $r_2, r_3, r_4, r_5, r_6$ est $\geqs 3$, compte-tenu du lien entre les $r_n$ et les $p_n$.

\

\noi Dans le m\me ordre d'id\Žes,  en 2003, {\sc Grekos, Haddad, Helou} et {\sc Pihko} [3] ont  \Žtabli ceci : pour toute base, et pour l'un au moins des nombres de repr\Žsentations associ\Žs, $r_n$, on a  $r_n\geqs6$. Ce r\Žsultat a \Žt\Ž am\Žlior\Ž, en 2006, par {\sc Borwein, Choi}, and {\sc Chu} [1], en $r_n\geqs 8$.

\

\head{Les polyn\™mes compliformes\\ Une parenth\se}

\

On a d\Žsign\Ž par M $= \Z[x_2,x_3,\dots]$ l'alg\bre de tous les polyn\™mes  aux variables $x_2,x_3,\dots$, \ˆ coefficients entiers. C'est l'alg\bre (commutative) libre engendr\Že par $x_2,x_3,\dots$, sur $\Z$.  On a d\Žsign\Ž par  C l'ensemble de tous les polyn\™mes {\it compliformes}, $u(x_2,x_3,\dots,x_n)$,  en $x_2,x_3,\dots$. C'est un sous-ensemble de l'alg\bre libre M. Ce n'est, en aucun cas, une sous-alg\bre ! C'est, cependant, un sous-module, un sous-groupe pour l'addition.

\

\noi D\Žtaillons un peu. \`A chaque polyn\™me $u(x_2,x_3,\dots,x_n)\in$ M correspond le polyn\™me compliforme $\hat u(x_2,x_3,\dots,x_n)$ obtenu en rempla\c cant  dans $u$ chaque $x_k^j$, pour $j\geqs 2$, par $x_k$. L'application $u\mapsto \hat u$ de M dans C est un homomorphisme sujectif de groupes additifs : c'est une {\bf r\Žtraction} de M sur son sous-groupe additif C. 
\su{Scholie} \sl Un  polyn\™me compliforme, $u(x,y,\dots,z)$, \ˆ $n$ variables est identiquement nul si, pour un choix de $n$ couples donn\Žs de valeurs  des variables, $x_0 \neq x_1, y_0\neq y_1, \dots, z_0 \neq z_1$, toutes les valeurs $u(x_i,y_j,\dots,z_k)$ du polyn\™me sont nulles\rm.

\

\noi {\bf En effet}, il va de soi, qu'un polyn\™me compliforme \ˆ une seule variable est identiquement nul si ses deux valeurs sont nulles  pour deux valeurs donn\Žes, distinctes, de la variable. De m\me, un  polyn\™me compliforme $u(x,y)$ \ˆ deux variables est identiquement nul si, pour des valeurs donn\Žes, $x_0 \neq x_1$ et $y_0\neq y_1$, on a $u(x_i,y_j) = 0$. Pour le voir, on \Žcrit 
$u(x,y) = v(y)x + w(y)$
o\ $u(y)$ et $w(y)$ sont compliformes. Puis on observe que les deux polyn\™mes compliformes \ˆ une variable, $u(x,y_0)$ et $u(x,y_1)$, sont  nuls, de sorte que $v(y_0), v(y_1), w(y_0), w(y_1)$ sont nuls donc $v(y)$ et $w(y)$ sont identiquement nuls. En poursuivant, ainsi, par r\Žcurrence, on arrive \ˆ la conclusion. {\bf cqfd}

\

\noi \sl En particulier, pour $u(x_2,x_3,\dots,x_n)\in$ M, les  \Žnonc\Žs suivants sont ainsi \Žquivalents :\rm

\

(i) \ \  \sl On a $u(x_2,x_3,\dots,x_n)= 0$ quels que soient les choix $x_k\in \{0,1\}$.\rm

(ii) \ \sl On a $\hat u(x_2,x_3,\dots,x_n)\equiv 0$.\rm

\

\noi En effet, pour $x_k\in \{0,1\}$, on a $u(x_2,x_3,\dots,x_n)=  \hat u(x_2,x_3,\dots,x_n)$, bien \Žvidemment. Le reste d\Žcoule de l'\Žnonc\Ž pr\Žc\Ždent, sans d\Žtour.

\

\noi On observera encore ceci. Un polyn\™me compliforme $u(x,y,\dots,z)$ est une fonction affine de chacune de ses variables;  on a ainsi :
$$u(s+t,y,\dots,z) = u(s,y,\dots,z)+ u(t,y,\dots,z)- u(0,y,\dots,z).$$

\

\head{L'alg\bre r\Žduite  de Erd\šs-Tur\'an\\
Un dernier avatar}

\

La forme de r\Žciprocit\Ž mise en \Žvidence ci-dessus ainsi que son corollaire montrent comment  on peut se passer, raisonnablement, des coefficients fractionnaires dans l'alg\bre ETB. Cela conduit \ˆ l'ultime simplification que voici.

\

\noi On appellera {\bf alg\bre r\Žduite de Erd\šs-Tur\'an}, et on d\Žsignera par ETR,  la sous-alg\bre r\Žduite de l'alg\bre ETB, engendr\Že  par $\Z$ et les $a_n$. C'est, tout simplement,  l'alg\bre de polyn\™mes $\Z[a_2,a_3,\dots]$ munie de la famille de relations d'idempotence : \[a_n^2 = a_n.\tag{$***$}\]
On observera que les $b_n$ et les $c_n$ appartiennent \ˆ cette alg\bre puisque l'on a 
\[b_n = a_0a_n + a_1a_{n-1} + \dots + a_ka_{n-k} + \dots +
a_na_0,\tag{$*$}\]
\[c_n =  \sum_{0\leqs k \leqs n/2} a_ka_{n-k}.\tag{5}\]
On aurait pu introduire l'alg\bre ETR plus t\™t, directement, sans les d\Žtours pr\Žc\Ždents. On a pr\Žf\Žr\Ž suivre le cheminement qui y a conduit !
On peut donner de nombreuses pr\Žsentations, diff\Žrentes, et \Žquivalentes, de  cette alg\bre ETR.  On peut la consid\Žrer  comme l'alg\bre   des polyn\™mes en $a_k,b_k,c_k,$ \ˆ coeffcicients entiers, $\Z[a_2,a_3,\dots;b_2,b_3,\dots;c_2,c_3,\dots]$,  munie des relations ($*$), ($**$) et ($***$), par exemple. Ou encore comme l'alg\bre de polyn\™mes $\Z[a_2,a_3,\dots;c_2,c_3,\dots]$ munie des relations  ($***$) et (5), dans laquelle  les \Žl\Žments $b_n$ seraient d\Žfinis par les formules ($*$). 

\

\noi Toujours est-il, l'alg\bre r\Žduite de Erd\šs-Tur\'an  est aussi le quotient de l'alg\bre libre M$=\Z[x_2,x_3,\dots]$ par son id\Žal $I$ engendr\Ž par les \Žl\Žments $x_n^2-x_n$. C'est, pour ainsi dire, une alg\bre de polyn\™mes \ˆ coefficients entiers, en des variables $a_2,a_3,\dots$ {\bf idempotentes}. On observe que tout \Žl\Žment non nul de l'id\Žal $I$ est un polyn\™me dont, au moins, l'une des variables figure avec un degr\Ž sup\Žrieur ou \Žgal \ˆ 2.

\

\head{Sur la structure de l'alg\bre ETR}

\

Bien entendu, les relations ($*$), ($**$), ($***$), (5) et (7), sont, toutes, valables dans l'alg\bre ETR.

\su{S1 Comment les $\mathbf {a_n}$ engendrent ETR} Le noyau de l'application suivante de C dans ETR :
$$u(x_2,x_3,\dots,x_n) \mapsto u(a_2,a_3,\dots,a_n)$$
 est nul car $u$ non nul  est de degr\Ž 1 au plus en chacune de ses variables. Cette application est donc un isomorphisme de groupes additifs et l'on a :
$$\mathrm{ETR} = \{u(a_2,a_3,\dots,a_n) : u\in \mathrm C\}.$$
Pour $u(x_2,\dots,x_n)\in$ C, on a, ainsi $u(a_2,a_3,\dots,a_n) = 0$ si et seulement $u(x_2,\dots,x_n) = 0$.

\

\noi Pour $u(x_2,x_3\dots x_n) \in$ M, on a  $u(a_2,a_3,\dots,a_n) = \hat u(a_2,a_3,\dots,a_n)$, bien entendu.

\

\noi On retiendra ceci. \`A chaque \Žl\Žment $z$ de l'alg\bre ETR est attach\Ž un unique polyn\™me compliforme $q_z = q_z(x_2,x_3,\dots,x_n)\in$ C, celui pour lequel on a $z = q_z(a_2,a_3,\dots)$. Lui est attach\Ž \Žgalement l'ensemble de tous les polyn\™mes $u(x_2,x_3,\dots,x_n)\in$ M pour lesquels on a $\hat u = q_z$. 

\su{S2 Universalit\Ž de ETR} Quel que soit l'anneau commutatif $R$ et la suite  d'idempotents, $i_2, i_3, \dots$, de  $R$, l'application $a_n\mapsto i_n$, pour $n\geqs 2$, se prolonge, de mani\re unique, en un homomorphisme de l'anneau ETR dans l'anneau $R$ : c'est  l'application $z= q_z(a_2,\dots,a_n) \mapsto q_z(i_2,\dots,i_n)$. On peut ainsi dire que l'alg\bre ETR est {\it universelle} dans la cat\Žgorie des anneaux commutatifs engendr\Žs par une suite (finie ou infinie) d'idempotents.

\

\su{S3 Absence d'int\Žgrit\Ž} L'alg\bre ETR n'est pas int\gre, loin de l\ˆ. Par exemple, on a  $a_n(a_n-1) = 0$,   alors que $a_n \neq 0$ et $a_n\neq 1$, pour tout $n\geqs 2$. On a aussi $(c_2 -1)(c_2-2) = 0$ tandis que $c_2\neq 1$ et $c_2\neq 2$. On a, plus g\Žn\Žralement,  $(c_{n+1} - e_n)(c_{n+1} - e_n-1) = 0$ alors qu'aucun des deux facteurs n'est nul.

\

\noi Voici encore une identit\Ž moins apparente. On a $c_3 = a_3 + a_2$ de sorte que, tous calculs faits, il vient $c_3(c_3-1)(c_3-2) = 0$.

\

\head{De la topologie en interm\de}

\

Parmi toutes les bases de $\N$, distinguons celles pour lesquelles la suite des entiers $r_n$ est  born\Že, et les autres : les premi\res sont celles qui ne satisfont pas la conjecture de Erd\šs-Tur\'an. S'il y en a, elles sont de {\it classe inf\Žrieure}, selon notre classification. On peut voir, sans grand effort, que si des bases de classe inf\Žrieure existent leur ensemble n'est pas bien grand, il est {\bf maigre} au sens de la th\Žorie de Baire. \{Comme on le sait, en topologie, la notion de partie maigre est analogue \ˆ celle de partie n\Žgligeable en th\Žorie de la mesure.\}

\

\noi Soit $S$ l'ensemble de toutes les s\Žries formelles de la forme 
$$u(t) = \sum u_nt^n \ \ \text{o\ $u_n \in \{0,1\}$ pour tout $n$}.$$ 
On identifie $S$ \ˆ l'espace compact $K = 2^\N$, en identifiant la suite $(u_n)_{n\in\N}$ au point correspondant dans $K$.  \`A chaque s\Žrie $u \in K$, on associe son carr\Ž $v(t) = u(t)^2 = \sum r_n t^n$. On d\Žsigne par $T$ le sous-espace de $K$ form\Ž des s\Žries $u$ pour lesquelles on a $r_n \geqs 1$ pour chaque indice $n$, et par $U$ le sous-espace de $T$ form\Ž des s\Žries $u \in T$ pour lesquelles la suite des $r_n$ est born\Že. Ainsi, en ce sens pr\Žcis, on peut assimiler  $T$ \ˆ l'espace des bases de $\N$ et $U$ au sous-espace des bases dont le profil est born\Ž. On montre, sans grand d\Žtour, que $U$ est une partie \sl maigre \rm dans l'espace $T$. 

\

\noi Cela veut dire que {\bf $\mathbf U$  est contenue dans une r\Žunion d\Žnombrable de parties ferm\Žes dont les int\Žrieurs sont vides dans $\mathbf T$.}

\

\su{D\Žmonstration} Soit $V = T \setminus U$. Pour chaque entier $p>0$, on pose
 $$V(p) = \{u \in T : \exists m  \ (r_m \geqs p)\}.$$
   Ainsi $V = \bigcap V(p)$. Or, d'une part, chaque $V(p)$ est ouvert dans $T$ : en effet, \Žtant donn\Ž $u \in V(p)$, on consid\re  $m$, l'un des indices pour lesquels on a $r_m \geqs p$; toute s\Žrie $h \in T$ dont les coefficients co\•ncident avec ceux de $u$ jusqu'au rang $m$ appartient encore \ˆ $V(p)$. D'autre part, la partie $V$ est dense dans $T$ : en effet, \Žtant  donn\Žs $h\in T$ et un indice $q$, on d\Žsigne par  $u$ la s\Žrie dont les $q$ premiers coefficients co\•ncident avec ceux de $h$, et dont les autres sont tous \Žgaux \ˆ $1$; cette s\Žrie $u$ appartient bien \ˆ $V$. Ainsi, $V$ est un \sl r\Žsiduel \rm dans $T$ et son compl\Žmentaire est donc \sl maigre\rm.\qed

\

\noi {\bf Dans l'ensemble des bases, celles qui v\Žrifient la conjecture de Erd\šs-Tur\'an sont \emph{g\Žn\Žriques}}, autrement dit, elle forment un r\Žsiduel.

\

\noi \{Cela suffit \ˆ expliquer pourquoi, si la conjecture de Erd\šs-Tur\'an \Žtait fausse, il est difficile de trouver des bases qui la contredisent. On pourra se souvenir de la situation analogue o\ l'on cherchait,  dans les ann\Žes 20 du si\cle dernier, des fonctions r\Želles continues qui ne poss\dent de d\Žriv\Že \ˆ gauche ou \ˆ droite (finie ou infinie), en aucun point; \ˆ ce propos, lire {\sc Saks} [10], ou voir l'{\sc Addendum} ci-dessous.\}

\head{Le syst\me principal}

Bien entendu, les relations ($*$), ($**$), ($***$), (5) et (7), sont, toutes, valables dans l'alg\bre ETR.

\

\noi  On reprend les relations (7) du lemme et, tout en gardant les m\mes notations [le risque de confusion demeurant minime], on les consid\re comme un syst\me d'\Žquations o\ les $a_n$ sont des donn\Žes et les $c_n$ sont  les inconnues. On \Žcrit ainsi :
\[c_{n+1} = a_{n+1} + e_n(c_2,\dots,c_n).\tag{SP}\]
En d\Žtaillant les premi\res \Žquations, il vient

\

$c_2 = a_2 + 1$,

$c_3 = a_3 + c_2 - 1$,

$c_4 = a_4 + c_3$,

$c_5 = a_5 + c_4 - 2c_3 - c_2+ c_3c_2 +1$,

$c_6 = a_6 +c_5 - 2c_4+ 4c_3+ c_4c_2 - 2c_3c_2$,

$c_7 = a_7  + c_6 - 2 c_5  + 4 c_4 - 4 c_3 + 3 c_2 +  c_5c_2 + c_4 c_3 - 3 c_4 c_2  + c_3 c_2  
 - 3$,
 
 $c_8= a_8 + c_7- 2 c_6+ 4 c_5 - 4 c_4 + 4 c_3  - 4 c_2+ c_6 c_2 + c_5 c_3- 3  c_5c_2- 2 c_4 c_3$
 
\qquad$+4 c_4 c_2 - c_3 c_2 +4 $,

$\dots\dots\dots\dots\dots\dots\dots$

$c_{n+1} = a_{n+1} + e_n(c_2,\dots,c_n)$,

$\dots\dots\dots\dots\dots\dots\dots$

\

\noi D'apr\s sa structure m\me, en sachant seulement que  les $e_n(c_2,\dots,c_n)$ sont des polyn\™mes compliformes, on peut d\Žduire de nombreuses propri\Žt\Žs du syst\me SP, dont voici les deux premi\res.

\su{P1 Existence des solutions} Quelles que soient les valeurs donn\Žes aux $a=(a_2, a_3, \dots, a_n\dots)$, le syst\me poss\de une solution unique  en $c[a]=(c_2, c_3, \dots, c_n\dots)$. Si les valeurs des $a_k$ sont des nombres r\Žels (respectivement, complexes), il en sera de m\me des valeurs des $c_k$. Si les valeurs des $a_k$ sont des entiers, les valeurs des $c_k$ seront des entiers. D'une mani\re g\Žn\Žrale, si les $a_k$ appartiennent \ˆ un anneau commutatif donn\Ž, quelconque, il en sera de m\me des $c_k$. 

\su{P2 Unicit\Ž du syst\me} Soit
\[c_{n+1} = a_{n+1} + t_n(c_2,\dots,c_n),\]
un autre syst\me
o\ les $t_n(c_2,\dots,c_n)$ sont encore des polyn\™mes compliformes. Si les deux syst\mes ont  les m\mes solutions, les polyn\™mes $t_n$ sont identiques aux polyn\™mes $e_n$. En fait, si l'on suppose (seulement)  que, pour chacun des choix  $a_n\in\{0,1\}$, les deux syst\mes ont la m\me solution $c[a]$, dans l'anneau $\Z$, alors les $t_n$ sont identiques aux $e_n$. Cela d\Žcoule, directement, de la scholie ci-dessus.

\

Au-del\ˆ de ses seules propri\Žt\Žs de structure, le syst\me SP poss\de des propri\Žt\Žs notables qui d\Žcoulent de la mani\re dont il a \Žt\Ž construit.

\su{P3 Les pr\Žsentations}   Soit $A$ une partie quelconque de $\N$ \ˆ laquelle appartiennent $0$ et $1$. On se replace dans l'alg\bre ETR et on prend $a_n = 1 \ou 0$ suivant que $n$ appartient ou non \ˆ $A$. La solution  en $c= (c_2,c_3,\dots)$ du syst\me SP qui correspond \ˆ la donn\Že $a= (a_2,a_3,\dots)$ est alors $c_n = p_n$, le nombre de $A$-pr\Žsentations de $n$. C'est comme cela que le syst\me a \Žt\Ž construit.

\su{Insistons}\`A chaque \Žtape, suivant que $n+1$ appartient ou non \ˆ $A$, on a
\[p_{n+1} = e_n(p_2,\dots,p_n) +1 \ou p_{n+1} = e_n(p_2,\dots,p_n).\tag{8}\]
Pour chaque indice, on a le choix entre deux valeurs diff\Žrentes de $p_{n+1}$ : l'entier $x_0 = e_n(p_2,\dots,p_n)\geqs 0$ et l'entier $x_1= x_0 + 1$. On peut imaginer que le choix se joue  {\it \ˆ pile ou face}. La suite des $p_n$ est ainsi engendr\Že d'une mani\re {\it presque} automatique. Le mot automatique ayant d\Žj\ˆ \Žt\Ž utilis\Ž en un sens pr\Žcis, pour d\Žcrire les suites engendr\Žes par un certain type d'automates, on utilisera (faute de mieux) un autre mot pour d\Žcrire les suites $1,p_2,p_3,\dots,$ associ\Žes aux parties $A\inc \N$ (auxquelles appartiennent $0$ et $1$) : on dira que ce sont les suites {\bf h\Žmitropes}. 

\

\noi Observons encore ceci. En particulier, $A$ est une base de $\N$ si et seulement si, pour tout $n$, on a $p_n\geqs1$. Pour faire court, on dira qu'une suite h\Žmitrope $1,p_2,p_3,\dots,$ est {\bf basique} lorsque l'on a $p_n \geqs 1$, pour tout $n\geqs1$.

\

\noi La conjecture de Erd\šs-Tur\'an revient   \ˆ dire qu'{\bf aucune suite h\Žmitrope basique n'est born\Že}.

\su{Une petite digression} Il existe des suites h\Žmitropes born\Žes, il en existe m\me qui sont major\Žes par $1$, et ne se terminent pas seulement par des $0$. Ces derni\res sont associ\Žes aux ensembles de Sidon.  En particulier, pour $A=\{0,1,3,\dots,\}= \{v_{n+1} = 2v_n + 1 : v_0 = 0, n\in \N\}$, la suite des $p_n$ est bien form\Že de $0$ et de $1$ et ne se termine pas seulement par des $0$.

\

\head{Un peu de mesure}

\

Marquons une pause. On peut montrer que, parmi les bases de $\N$, il est infiniment peu probable qu'il y en ait une qui contredise la conjecture de Erd\šs-Tur\'an. Dans une correspondance priv\Že,  {\sc Jean-Pierre KAHANE} m'a montr\Ž avec pr\Žcision comment on \Žtablit ce genre de r\Žsultat.

\

\noi On se donne un nombre r\Žel $p\in \ ]0,1]$ et on suppose que $P(a_n=1) = p$ et $P(a_n = 0) = 1-p = q$. Cela veut dire que, dans un tirage au sort pour d\Žterminer $a_n\in \{0,1\}$, on suppose que la probabilit\Ž de tirer le $1$ est \Žgale \ˆ $p$. Cela revient \Žgalement \ˆ dire que l'on munit l'ensemble discret $\{0,1\}$ de la mesure de probabilit\Ž $P$ pour laquelle $P(\{1\}) = p$ et $P(\{0\})= q$ puis l'on munit l'espace-produit $K = \{0,1\}^\N = 2^\N$ de la mesure-produit. On introduit les variables al\Žatoires suivantes :
$$c_n  = \sum_{0\leqs i\leqs n/2} a_ia_{n-i}= a_0a_n + a_1a_{n-1} + \cdots.$$
Pour chaque entier $k$, on montre alors que 
$$P\left(\sup_{n\geqs 0} c_n \leqs k\right) =0,$$
ce qui revient \ˆ dire que, {\bf presque s\žrement}, la suite des $c_n$ est non born\Že. L'ensemble des bases de $\N$ qui contredisent la conjecture de Erd\šs-Tur\'an est donc {\bf n\Žgligeable}, quel que soit le choix d'une probabilit\Ž $p>0$.
\su{D\Žmonstration}  L'\Žv\Žnement $c_{2m+1} = k$ se produit si et seulement s'il existe une partie $I\inc \{0,1,\dots,m\}$ ayant $k$ \Žl\Žments et pour laquelle on a
$$a_ia_{2m+1-i} = 1, \pour i\in I, \et a_ia_{2m+1-i} = 0,\sinon\nts.$$
Cela revient \ˆ dire que l'ensemble $I = \{i : a_i = a_{2m+1-i} = 1 \et i\leqs m\}$ est form\Ž de $k$ \Žl\Žments. On a ainsi :
$$P(c_{2m+1} = k) = \binom{m+1} k p^{2k}(1-p^2)^{m+1-k}.$$
De m\me, l'\Žv\Žnement $c_{2m} = k$ se produit si et seulement s'il existe une partie $I \inc \{0,1,\dots,m\}$ ayant $k$ \Žl\Žments et pour laquelle on a
$$a_ia_{2m-i} = 1, \pour i\in I, \et a_ia_{2m-i} = 0,\sinon\nts.$$
En tenant compte des deux cas possibles, selon que $m$ appartient ou non \ˆ la partie $I$, on a
$$P(c_{2m} = k) = \binom m {k-1} p^{2k-1}(1-p^2)^{m+1-k} + q\binom m k p^{2k}(1-p^2)^{m-k}.$$
On v\Žrifie simplement, \ˆ l'aide de la formule de Stirling, que $P(c_n=k)$ tend vers $0$, pour $k$ fix\Ž, lorsque $n$ tend vers l'infini. D'o\ le r\Žsultat.\qed

\

\noi Dans le m\me ordre d'id\Žes, et en poussant un peu plus loin les calculs, on peut montrer ceci : quel que soit le choix que l'on fait d'un $p> 1/2$, on aura
$$\limsup_{n} \frac{8c_n}{n}\geqs 1, \ \text{presque s\žrement}.$$

\

\head{La pause}

\

On a vu, plus haut, comment et pourquoi la topologie laisse tr\s peu de place \ˆ d'\Žventuelles bases qui contrediraient la conjecture de Erd\šs-Tur\'an. On vient de voir que la probabilit\Ž (le hasard) ne leur accorde aucune chance. Ainsi, la recherche d'un contre-exemple \ˆ la conjecture {\bf\emph{ET}} risquerait-elle fort de ressembler \ˆ  {\sl une marche sur des \oe ufs}, \ˆ un passage \ˆ travers de nombreuses {\sl portes \Žtroites}, \ˆ un cheminement sur {\sl une route de cr\te} bord\Že de pr\Žcipices. 

\

\noi Que peut encore nous apprendre l'alg\bre ? ou la combinatoire ?

\

\head{L'arbre binaire des suites h\Žmitropes} 

\

 L'ensemble des  suites h\Žmitropes peut \tre consid\Žr\Ž comme un arbre binaire $\cal T$. Les n\oe uds de l'arbre $\cal T$ sont les suites finies $(1,p_2,\dots,p_n)$  o\ $(1,p_2,p_3,\dots)$ est une suite h\Žmitrope. La racine de cet arbre est la suite   \ˆ un seul terme $(p_1=1)$. Chaque n\oe ud $(1,p_2,\dots,p_n)$ a exactement deux successeurs imm\Ždiats,  un  {\bf n\oe ud inf\Žrieur}, $(1,p_2,\dots,p_n,x)$, 
et   un {\bf n\oe ud sup\Žrieur} $(1,p_2,\dots,p_n,x+1)$. Les branches (infinies) de l'arbre sont les suites h\Žmitropes
 $(1,p_2,\dots,p_n,\dots)$.
 
 \
 
 \noi Un n\oe ud  sera dit {\bf persistant} lorsqu'aucun de ses $p_k$ n'est un $0$. De m\me, une branche sera dite {\bf persistante}  lorsqu'aucun de ses $p_k$ n'est un $0$. Dans l'arbre $\cal T$, il y a toujours, au moins, une branche persistante passant par un n\oe ud persistant donn\Ž, quelconque :  en effet, apr\s chaque n\oe ud, il suffit de choisir le successeur imm\Ždiat sup\Žrieur pour obtenir une branche persistante.

 \

\noi \`A chaque
branche $p=(1,p_2, \dots)$ correspond la partie $A \inc \N$ pour laquelle $p_n$ est le nombre de $A$-pr\Žsentations de $n$. La branche $p$ est persitante si et seulement si $A$ est une base de  $\N$ : cela \Žtablit une relation  bijective  entre  les branches persistantes et les bases.

\

\noi La conjecture {\bf \emph{ET}} s'\Žnonce alors comme suit :
L'arbre $\cal T$ ne poss\de pas de branche persistante  $p=(1,p_2, \dots)$ pour laquelle l'ensemble des $p_n$ est born\Ž.

\

\noi Il y a plusieurs mani\res diff\Žrentes pour {\it coder} les n\oe uds et les branches de notre arbre. En voici deux partiuli\res : l'une utilise des $0$ et des $1$, l'autre utilise des $+$ et des $-$.

\su{Le code en $\mathbf 0 \et \mathbf 1$} Le code du n\oe ud $(1,p_2,\dots,p_n)$ s'obtient en rempla\c cant $p_k$ par $0$ si $p_k$ est pair, par $1$ si $p_k$ est impair.
\su{Le code en $\mathbf + \et \mathbf -$} On obtient le code du n\oe ud $(1,p_2,\dots,p_n)$  en rempla\c cant $p_k$ par $+$ ou $-$ suivante que 
$(1,p_2,\dots,p_k)$ est sup\Žrieur ou inf\Žrieur.

\

\noi Chacun de ces deux codes est univoque, bien entendu.

\

\head{Addendum} 

\

Dans les ann\Žes 20 du si\cle dernier, au sujet de l'existence de fonctions continues sans d\Žriv\Že, on envisageait deux classes de  fonctions, en particulier.   La premi\re classe est celle des fonctions r\Želles continues d\Žfinies sur l'intervalle $I = [0,1]$ qui, en chaque point, n'ont pas de d\Žriv\Že  \ˆ droite finie, ni de d\Žriv\Že \ˆ gauche finie. Dans l'espace de Banach $\cal C(I)$, on savait montrer que l'ensemble de ces fonctions avait  {\bf un compl\Žmentaire  maigre} [{\sc Banach} et {\sc Mazurkiewicz}]. On pouvait donc dire qu'il en existait beaucoup sans avoir besoin d'en montrer une seule. La seconde classe est celle des fonctions r\Želles continues d\Žfinies sur l'intervalle $I = [0,1]$ qui, en chaque point, n'ont pas de d\Žriv\Že  \ˆ droite {\bf finie ou infinie}, ni de d\Žriv\Že \ˆ gauche {\bf finie ou infinie}. {\sc Saks} a montr\Ž que l'ensemble de ces fonctions \Žtait {\bf une partie maigre} de l'espace de Banach $\cal C(I)$ expliquant ainsi la raison pour laquelle on  avait mis beaucoup de temps   avant [{\sc Besicovitch}] d'en exhiber une !

\

Voici les r\Žf\Žrences.

\

\noi Banach, S. Ueber die Bairesche Kategorie gewisser Funktionenmengen. Stud. Math. 3 (1931) 174-179.

\

\noi Mazurkiewicz, S.
Sur les fonctions non d\Žrivables.  Stud. Math. 3 (1931) 92-94.

\

\noi Besicovitch, A. Diskussion der stetigen Funktionen in Zusammenhang mit der Frage \Ÿber ihre Differentiierbarkeit.II. Bull. de l'Ac. des Sc. de Russie, 19 (1925) 527-540.

\

\head{En guise de conclusion provisoire}

\

En un sens bien pr\Žcis, on a vu qu'il y avait peu de place et presque aucune chance pour que la conjecture de Erd\šs-Tur\'an soit fausse. Peu de place ? Ou pas du tout !

\centerline{\bf \`A suivre}

\newpage

\centerline{\bf \small Voici la version en lange anglaise}

\

\head{\bf \Large Some peculiarities of order 2 bases of $\N$\\ and the Erd\šs-Tur\'an conjecture}

\

\hfill{\sl If you had known the virtue of the ring,} 

\hfill{\sl Or half her worthiness that gave the ring,} 

\hfill {\sl Or your own honour to contain the ring,} 

\hfill{\sl You would not have parted with the ring...}

\hfill{The Merchant of Venice,
V, I.199-202}

\

\head{Introduction}

\

\

Given a subset $X$ of $\Z$ and an integer $z$, let $p(X,z)$ be the number of all ordered couples $(x,y)\in X \times X$  such that $x+y=z$\nts: Of course, this can either be an integer or infinity.

\

\noi In 2005, among many other more general results, {\sc Nathanson} [9] showed the following\nts: Let $(s_z)_{z\in \Z}$ be a (double) sequence of integers $s_z>0$. Then, there always exists a subset $X\inc \Z$ such that $p(X,z) = s_z$, for each $z$.

\

\noi Now, take  a subset $A\inc\N= \{0,1,2,\dots\}$ and, for the sake of brevity, use the following terminology and notation\nts:  Call each ordered couple $(x,y)\in A\times A$ such that $x+y = n$ an $A$-{\it representation} of $n$, and let $r_n$ be the number of all those $A$-representations. The subset $A$ is called a {\it basis} for $\N$ whenever $r_n > 0$ for each $n\in \N$. When this is the case, it is easily seen that $0$ and $1$ must belong to $A$, so that one has\nts:
$$r_0 = 1, r_1 = 2.$$
It can also as easily be seen that
$$r_2 = 1 \ora 3 \ , \ r_3 =r_2- 1 \ora r_2 + 1 \ , \ r_4 =  \frac{2r_3-r_2+1}{2}  \ora  \frac{2r_3- r_2+3}{2}.$$
As we shall presently show, more generally, there exists a sequence  of polynomials, $d_n(b_2,\dots,b_n)$, such that, for each basis $A$ and each $n\geqs 1$, one has\nts: 
$$r_{n+1} = d_n(r_2,\dots,r_n) \sia n+1\notin A,$$
$$r_{n+1} = d_n(r_2,\dots,r_n) + 2 \sia n+1\in A.$$
Their degrees  in each of the variables is equal to $1$, so that they are multilinear functions of the variables, and their coefficients are in the ring $\Z[1/2]$. For instance, we have\nts:
$$d_1 = 1, \ d_2(b _2) = b _2-1, \ 2d_3(b _2,b _3) = 2b _3 - b _2+1,$$
$$2d_4(b _2,b _3,b _4) = b _4  -\ 3b _3 - b _2 + b_3b_2 +  1,$$
$$2d_5(b _2,b _3,b _4,b _5) = 2b _5 - 3 b _4  + 5 b _3 + b _2+ b _4b _2-2  b _3b _2  - 1,$$
$$4d_6 = 4b _6 - 6 b _5 + 10 b _4 - 13 b _3+ 6 b _2 + 2b _5b _2  + 2 b _4 b _3 - 6 b _4 b _2 + 3 b _3 b _2   - 6.$$
Clearly, this result is in contrast with {\sc Nathanson}'s. The sequence $(r_n)_{n\in \N}$, far from being arbitrary, is much {\bf constrained}. Let us insist\nts: This sequence does not unfold freely, {\it it is subject to high monitoring}.

\

\noi In 1941, {\sc Erd\šs} and {\sc Tur\'an} [2] made a conjecture which can be stated as follows\nts:
\su{The \emph{ET}  conjecture}{\sl For each basis $A$ of $\N$, the sequence $(r_n)_{n\in \N}$  of the $A$-representation numbers is unbounded.}

\

\noi In the studies about this conjetcure, a well known method is frequently used, (see, for instance, [1], [3], [4], [6], [7], [8], [11])\nts: To each subset $A\inc\N$, two formal series, $f(t)$ and $g(t)$, are associated. The first series represents $A$ as \nts:
$$f(t) = \sum_{a\in A} t^a.$$
It takes the following form\nts:
$$f(t) = \sum_{n\geqs 0} u_nt^n$$
where $u_n = 1 \ora  0$, according as  $n$ belongs or not to $A$. The second series is the square of the first, $g(t) = f(t)^2$. It is easily seen that the following holds\nts:
$$g(t) = \sum_{n\geqs 0} r_nt^n$$
where  the $r_n$'s are precisely the numbers defined above. When $A$ is a basis, one must have $u_0 = u_1 = 1$, that is, $0$ and $1$ must belong to $A$, therefore, $r_0 = 1$ and $r_1=2$, as we already observed above.

\

\noi Define an $A$-{\it presentation} of $n$ to be any ordered couple $(x,y)\in A \times A$ such that $x+y=n$ and $x\leqs y$, then let $p_n$ denote the number of all those $A$-presentations. We shall call the sequence $p=(p_0,p_1,p_2,\dots,p_n,\dots)$ {\bf the profile} of the subset $A$, in order to avoid periphrases. Although close, and related, $A$-presentations and $A$-representations are not the same. In fact, the easy relation is the following\nts: If $n/2$ belongs to $A$, then $r_n = 2p_n - 1$; otherwise, $r_n = 2p_n$. So, when $n$ is odd, we always have $r_n = 2p_n$.  Sometimes, it seems wiser to use the $p_n$'s rather than the $r_n$'s, as we shall see, later.

\

\noi This study, quite naturally, leads to  the introduction of an algebra which I call {\it the Erd\šs-Tur\'an algebra} (or ring) and denote ET. Its definition follows.

\

\head{The Erd\šs-Tur\'an algebra}

\

Let $\mathrm D = \Z[1/2] = \{k/2^n : k\in \Z, n\in \N\}$ be the ring of dyadic rational numbers. Three sequences of variables are introduced, $(a_1, a_2, \dots)$, $(b_1, b_2, \dots)$, $(c_1, c_2,\dots)$. Let L be the ring of polynomials in those variables with coefficients in D, i.e.,
$$\mathrm{L = D}[a_1, a_2, \dots, b_1, b_2, \dots, c_1, c_2,\dots].$$
This is the free algebra generated on D by the variables. The following three formal series in $t$ with coefficients in L are introduced\nts:
$$f(t) = 1 + a_1t + a_2t^2 + \dots = 1 + \sum_{n\geqs 1} a_nt^n,$$
$$g(t) = 1 + b_1t + b_2t^2 +\dots = 1 + \sum_{n\geqs 1} b_nt^n,$$
$$h(t) = 1 + c_1t + c_2t^2 +\dots = 1 + \sum_{n\geqs 1} c_nt^n.$$
Forcing the  following two identifications
\[g(t) = f(t)^2 \anda 2h(t) = g(t) + f(t^2)\tag{0}\]
obtains an infinite system of relations among the variables. More precisely, we thus have two types of relations. Setting $a_0 = 1$ and computing the square of the series $f(t)$, we get the first type of relations\nts:
\[b_n = a_0a_n + a_1a_{n-1} + \dots + a_ka_{n-k} + \dots +
a_na_0.\tag{$*$}\]
Relations of the second type come from the second identification, the link between $f$, $g$, and $h$. They read\nts:
\[2c_{2n+1} = b_{2n+1} \anda 2c_{2n} = b_{2n} + a_n.\tag{$**$}\]
The  algebra ET is, then, defined to be the algebra $L$ subject to the whole set  of relations ($*$) and ($**$).

\ 

\noi In the formulas ($**$), both cases, even and odd, can be recast into one. Just agree to set $a_x = 0$, for each $x\notin \N$, to get :
\[2c_n = b_n + a_{n/2}.\tag{$**$}\]
Formulas ($*$) are simple enough and determine the $b_n$'s as functions of the $a_k$'s. Similarly, there are formulas (not quite as simple) which
determine the 
$a_n$'s as functions of the $b_k$'s. Somehow {\it more hidden}, those
formulas come as a blend of {\bf Catalan numbers} and
ordinary {\bf Bell polynomials}. In order to obtain them, we make the most of the square roots of formal series. Let us elaborate.

\

\head{Square roots of formal series}

\

\su{A1 Catalan numbers} The Catalan numbers
$C_1,C_2,\dots,C_k,\dots,$ are defined as follows\nts:
$$C_1 = 1\ \et \ C_k = C_1C_{k-1} + C_2C_{k-2} +
\dots C_{k-1}C_1, \pour k \geq 2.$$ 
So, they are {\bf positive integers}.  The first few values are\nts:
$$C_1 = 1 \ , \ C_2 = 1 \ , \ C_3 = 2 \ , \ C_4 = 5 \ , \
C_5 = 14.$$
Due to the recurrence relation, their generating sequence, the formal series $C(t) = \sum_{k \geqs 1} C_kt^k$,
satisfies the identity  $C = t + C^2$.  We thus have
$$C(t) = \frac{1 - \sqrt{1 - 4t}}{2} \ , \ \text{so that}
\
\ \sqrt{1 - 4t} = 1
 - 2\sum_{k \geqs 1} C_k t^k.$$ 
By the binomial formula, we get
$$ \sqrt{1 - 4t} = 1 + \sum_{k
\geqs 1}
\binom {1/2}{k}(-4)^kt^k.$$
Whence,
$$C_k = \frac{(-1)^{k-1}}{2} 4^k \binom {1/2}{k} =
\frac{1}{k}\binom {2k - 2}{k -1}.$$ 
So, by Wallis' formula, we  have$$\lim_{n \to \infty} 4 \frac{C_k}{4^k} k^{3/2} =
\frac{1}{\sqrt\pi}.$$

\su{A2  Bell polynomials} The
(ordinary) Bell polynomials $P_{n,k}$,  are
given by the defining relations 
\[(b_1t + \dots + b_nt^n + \cdots)^k =  P_{k,k}t^k +
P_{k+1,k}t^{k+1} + \dots + P_{n,k}t^n +
\cdots,\tag {1}\]
for integers $n \geqs k
\geqs 0$. In all other cases, we set
$P_{n,k} =0$, as a convention. 

\

\noi  It should be noticed that all the coefficients
of the polynomial $P_{n,k}$, in the variables $b_1,b_2,\dots,b_{n-k+1}$ , are positive integers as is clearly seen from the defining relations.
We shall abbreviate 
$P_{n,k}(b_1,b_2,\dots,b_{n-k+1})$ into
$P_{n,k}[b]$,  for  simplicity sake, it being assumed  that we  set $b = b_1t + \dots + b_nt^n + \cdots$.

\

\noi Remember that,  by definition, we have set $P_{n,k} = 0$ when the inequalities $n \geqs k \geqs 0$ do not hold\! ! This brings great simplifications in the handling of
summation indices.

\su{A3 Square roots of formal series} Consider again the formal series $f(t)$ with coefficient in the algebra $ET$\nts:
$$f(t) = 1 + a_1t + \dots + a_nt^n +
\cdots,$$  
and its square
$$g(t) = f^2(t) = 1 + b_1t + \dots + b_nt^n + \cdots.$$
We get  the simple formulas\nts:
\[b_n = a_0a_n + a_1a_{n-1} + \dots + a_ka_{n-k} + \dots +
a_na_0.\]
Here are the (reciprocal) formulas that determine the 
$a_n$'s as  functions of the  $b_k$'s : 
\[a_n = 2 \sum_{1 \leqs k \leqs n}
\frac{(-1)^{k-1}}{4^k} C_k
P_{n,k}(b_1,b_2,\dots,b_{n-k+1}).\tag{2}\]

\su{Proof} Simply, set $b = b(t) = b_1t +
\dots + b_nt^n + \cdots,$ and use the binomial formula to get\nts:
$$f(t) = (1+b)^{1/2} = 1 + \sum_{k \geqs 1} \binom{1/2}{k}
b^k.$$ 
Then, \sl \ˆ  la Faa di Bruno\rm, substitute to
$b^k$ its expansion in (1), getting
$$1 + \sum_{n \geq 1} a_nt^n = 1 + 2 \sum_{n \geq 1} \sum_{k
= 1}^n
\frac{(-1)^{k-1}}{4^k}C_k P_{n,k}[b]t^n.$$
Formulas (2) are obtained by identifying the coefficients of
$t^n$ on both sides.\qed

\

\noi For more details,  see, for example,  {\sc Louis COMTET}, 
Advanced combinatorics. Transl. from the French by J. W. Nienhuys. (1974)

\su{A4 Multi-indices} We shall use {\bf multi-indices}, very helpful when dealing with Bell polynomials. A multi-index, more precisely, a {\bf $k$-index} is 
a $k$-sequence  (a finite sequence of $k$ terms), $\mu =
(\mu_1,\dots,\mu_k)$ of integers $\mu_j > 0$.
To a multi-index we associate its {\bf weight}, the integer
$|\mu| = \mu_1 + \mu_2 +\dots + \mu_k$. Given a series
$b = b_1t + b_2t^2 + \dots + b_nt^n + \cdots$, we also associate to a $k$-index
$\mu = (\mu_1,\dots,\mu_k)$ the {\bf monomial}
$b_\mu = b_{\mu_1}b_{\mu_2}\dots b_{\mu_k}$. Let  $M(k,n)$ be the set of all $k$-indices having a weight equal to
$n$. Then, as is easily seen from their defining relations, Bell polynomials can be written  as\nts:  
$$P_{n,k}[b] = \sum_{ \mu \in M(k,n)}
b_\mu.$$
For instance, we therefore have\nts:
$$P_{n,1}[b] = b_n \ , \  P_{n,2}[b] = b_1b_{n-1} + b_2
b_{n-2} + \dots + b_{n-1} b_1 \ , \   P_{n,n}[b] = b_1^n.$$
By ($*$), when the $a_n$'s are integers, so are the
$b_n$'s, therefore, equally well, the $b_\mu$'s.

\su{A5 Coefficients of the square root series} Formula (2) can be written as\nts:
$$2a_n = 2a_n[b] = 4\sum_{\mu \in M(k,n)}
\frac{(-1)^{k-1}}{4^k} C_k  b_\mu.$$  
For instance, we have\nts:

\

$2a_1  = b_1$,

$2a_2  = b_2 - \frac{1}{4}b_1^2$,

$2a_3  = b_3 - \frac{1}{2} b_1b_2 +
\frac{1}{8}b_1^3$,

$2a_4 = b_4 - \frac{1}{2}b_1b_3 - \frac{1}{4}b_2^2 + \frac{3}{8}b_1^2b_2 - \frac{5}{64}b_1^4$,

$2a_5 = b_5 -\frac{1}{2}b_1b_4 - \frac{1}{2}b_2b_3 + \frac{3}{8}b_1^2b_3 + \frac{3}{8}b_1b_2^2 - \frac{5}{16}b_1^3b_2 + \frac{7}{128}b_1^5$,

$2a_6 = b_6 - \frac{1}{2}b_1b_5 - \frac{1}{2}b_2b_4 + \frac{3}{8}b_1^2b_4 - \frac{5}{16}b_1^3b_3- \frac{1}{4}b_3^2 -\frac{15}{32}b_1^2b_2^2  +  \frac{3}{4}b_1b_2b_3 + \frac{1}{8}b_2^3$ 

\qquad \ \ $+ \frac{35}{128}b_1^4b_2- \frac{21}{512}b_1^6$,

\

$\dots\dots\dots\dots\dots\dots\dots$

\[2a_n  = b_n - \frac{1}{4}(b_1b_{n-1} + b_2
b_{n-2} + \dots + b_{n-1} b_1)
+ 4\sum_{\underset{ k > 2}{\mu\in M(k,n)}}
\frac{(-1)^{k-1}}{4^k} C_k  b_\mu.\tag{3}\]

\

\noi This produces the following noteworthy relation on which we shall  comment, later\nts:
\[b_n = 2a_n + \frac{1}{4}(b_1b_{n-1} + b_2
b_{n-2} + \dots + b_{n-1} b_1)
+ 4\sum_{\underset{ k > 2}{\mu\in M(k,n)}} \frac{(-1)^k}{4^k}C_k  b_\mu.\tag{4}\]

\

\noi Observe that the term
$\sum_{\underset{ k > 2}{\mu\in M(k,n)}}
\frac{(-1)^{k-1}}{4^k} C_k  b_\mu$
is a polynomial depending on $b_1,\dots,b_{n-2}$ but not $b_{n-1}$\nts: Indeed, for  $k>2$ et and each of the $k$-indices  $\mu= (\mu_1,\dots,\mu_k)$ having weight equal to $n$,  one has  $\mu_j \leqs n-2$, by necessity. Let us write\nts:
$$\sum_{\underset{ k > 2}{\mu\in M(k,n)}}
\frac{(-1)^{k-1}}{4^k}C_k  b_\mu = Q_n(b_1,\dots,b_{n-2}).$$  

\su{\emph{En passant}} We notice, on the way, that the  Erd\šs-Tur\`an  conjecture {\bf\emph{ET}} can, now, be stated as follows\nts : For each given bounded sequence of integers  $b_k\geqs 1$, there is an index $n$ such that the $a_n$ in  formula (3)  is not a $0$ nor a $1$.

\su{A6 A special case} \sl The case where each of the $a_n$'s is equal to $0$ or $1$\rm.

\noi In that case, we have $a_n^2 = a_n$, and this relation implies the following\nts:
$$\left(\frac{1}{2}b_n - \frac{1}{8}(b_1b_{n-1} + b_2
b_{n-2} + \dots + b_{n-1} b_1)
+ 2Q_n(b_1,\dots,b_{n-2})\right)^2 =$$
$$\frac{1}{2}b_n - \frac{1}{8}(b_1b_{n-1} + b_2
b_{n-2} + \dots + b_{n-1} b_1)
+ 2Q_n(b_1,\dots,b_{n-2}).$$ 
An expression of $b_n^2$ obtains as a polynomial in  $b_1, b_2, \dots, b_n$, having degree $1$ in $b_n$.
Starting from
$$b_{n+1} = 2a_{n+1} + \frac{1}{4}(b_1b_{n} + b_2
b_{n-1} + \dots + b_{n} b_1) + 4 Q_{n+1}(b_1,\dots,b_{n-2},b_{n-1}),$$
substituting to $b_{n-1}^2, b_{n-2}^2,\dots$, {\it in order},  their respective values, we get an expression of $b_{n+1}$ of the form
$$b_{n+1} = 2a_{n+1} + d_n$$
where $d_n$ is a polynomial in $b_1, b_2, \dots, b_n,$ each variable having degree $1$, at most.

\

\noi We are led to a  {\it specialisation} of  the algebra ET as follows.

\

\

\head{The special Erd\šs-Tur\'an algebra}

\

I call {\bf special Erd\šs-Tur\'an algebra}, and denote ETS, the algebra obtained, starting from ET,   by adding the following sequence of idempotence relations\nts: 
\[a_n^2 = a_n \ , \ n = 1, 2, \cdots.\tag{$***$}\]
That is, ETS is the free algebra L with the  relations ($*$), ($**$), and ($***$).

\

\noi Of course, from the fundamental relations, many others are derived.

\

\noi For instance, putting together ($**$) and ($*$), obtains the following identities\nts:
$$c_{2n} = a_{2n} +  a_{2n-1} + a_2a_{2n-2} + \dots + a_{n-1}a_{n+1} + a_n,$$
$$c_{2n+1} =  a_{2n+1} +  a_{2n} + a_2a_{2n-1} + a_3a_{2n-2}+\dots + a_na_{n+1}.$$
Setting $a_0=1$, they both are held in the following\nts:
\[c_n = \sum_{0\leqs k\leqs n/2} a_ka_{n-k}.\tag{5}\]

\noi Still more, mixing (3) and ($***$), we get an expression of the form
\[b_{n+1} = 2a_{n+1} + d_n(b_1,b_2,\dots,b_n)\tag{6}\]
where $d_n$ is a polynomial in the variables $b_1, b_2, \dots, b_n,$ each variable having degree $1$, at most. As in the special case presented above, there is a (simple) algorithm  which  computes  each of the polynomials $d_n$.

\su{A general reduction algorithm} This is an {\bf elimination} procedure. Let $x_1,x_2,\dots$, be a sequence of variables given with the following relations\nts: 
$$x_n^2 = t_n(x_1,\dots,x_n),$$
where $t_n(x_1,\dots,x_n)$ is a  polynomial in $x_1,\dots,x_n,$ having degree  one, at most, in the variable $x_n$. Each polynomial $P(x_1,\dots,x_n)$ can be reduced to a polynomial in $x_1,\dots,x_n,$ with each variable having degree one, at most, by elimination of the powers $x_k^j$, for $j\geqs 2$.

\

\noi Indeed, for each integer $j\geqs 2$, the monomial $x_n^j$ has an expression as a polynomial in $x_1,\dots,x_n,$ having degree one in $x_n$. This can be done by induction, starting with $x_n^3 = t_n(x_1,\dots,x_n)x_n$. 

\

\noi Substituting, in the polynomial $P(x_1,\dots,x_n)$, to each of the monomials $x_n^j$ its expression obtains a polynomial $P_n(x_1,\dots,x_n)$ in which the variable $x_n$ appears with degree one, at most. Repeating the process gives the polynomials  $P_{n, n-1}, P_{n, n-1, n-2}, \dots, P_{n, n-1,\dots,1}$, and the result.

\su{Another alternative} In the preceding algorithm, instead of substitutions, one can, simply, use  Euclidean division. Specifically, divide the polynomial $P(x_1,\dots,x_n)$ by $T(x_n) = x_n^2 - t_n(x_1,\dots,x_n)$ for the variable $x_n$\nts :
$$P(x_n) = Q(x_n)T(x_n) + R(x_n).$$
The remainder thus obtained, $R(x_n)$, is precisely the polynomial $P_n$. This can be done at each step $n-1,n-2,\dots,1$.

\

\su{Algebras in a tower} We thus have a tiered structure of polynomial algebras\nts: 
$$\mathrm L = \mathrm D[a_1, a_2, \dots, b_1, b_2, \dots, c_1, c_2,\dots]$$
$$\mathrm{ET   = L}\  \  \mathrm{modulo} \ \ (*) \et (**)$$
$$\mathrm{ETS = ET} \  \ \mathrm{modulo} \ \ (***)$$

We, now, add one more level, {\it for the bases}, specializing as follows\nts:
$$\mathrm{ETB = ETS}\  \ \mathrm{modulo}\ \  a_1 = 1.$$

\

\

\head{The algebra ETB\\
for bases}

\

In the case of bases, we always have $r_1 = 2$. We, therefore, specialize the algebra ETS into ETB, by setting $a_1 = 1$, $b_1=2$\nts: This new algebra is a tool  {\bf well suited} to the study of the numbers $r_n$  and $p_n$ of representations and  presentations associated to bases.  

\

\centerline{\{In what follows, we always set $a_0 = a_1 = c _0 = c_1= 1$\}}

\

\noi In that algebra, there are many interesting derived  relations; here are some of the first; they are special cases of (3)\nts:

\

$2a_2  = b_2 - 1$,

$2a_3  = b_3 -  b_2 + 1$,

$2a_4 = b_4 - b_3 - \frac{1}{4}b_2^2 + \frac{3}{2}b_2 - \frac{5}{4}$,

$2a_5 = b_5 -b_4 - \frac{1}{2}b_2b_3 + \frac{3}{2}b_3 + \frac{3}{4}b_2^2 - \frac{5}{2}b_2 + \frac{7}{4}$,

$2a_6 = b_6 - b_5 - \frac{1}{2}b_2b_4 + \frac{3}{2}b_4 - \frac{5}{2}b_3+ \frac{35}{8}b_2 - \frac{1}{4}b_3^2 -\frac{15}{8}b_2^2  +  \frac{3}{2}b_2b_3+ \frac{1}{8}b_2^3 - \frac{21}{8}$.

\

\noi Write $a_2^2 = a_2$ and obtain $b_2^2 = 4b_2 -3$. Substitute to $b_2^2$ its value in the expressions of $2a_4$ and $2a_5$, to get

\

 $2a_4 = b_4 - b_3 - \frac{1}{4}(4b_2-3) + \frac{3}{2}b_2 - \frac{5}{4}= b_4 - b_3+\frac{1}{2}b_2 -\frac{1}{2}$.
 
 $2a_5 = b_5 -b_4 - \frac{1}{2}b_2b_3 + \frac{3}{2}b_3 + \frac{3}{4}(4b_2 - 3)- \frac{5}{2}b_2 + \frac{7}{4}= b_5 -b_4 -\frac{1}{2}b_2b_3+$
 
  \qquad \ \ $\frac{3}{2}b_3 +  \frac{1}{2}b_2- \frac{1}{2} $.

\

\noi Write $a_3^2 = a_3$ and obtain $b_3^2 = 2b_2b_3 -4b_2 + 4$, $b_2^3 = 13b_2 - 12$. Substitute to $b_2^2$, $b_2^3$, and $b_3^2$ their respective values in the expression of $2a_6$, to get

\

$2a_6 = b_6 - b_5 + \frac{3}{2}b_4 - \frac{1}{2}b_2b_4  - \frac{5}{2}b_3+ \frac{35}{8}b_2 - \frac{1}{4}(2b_2b_3 -4b_2 + 4)-\frac{15}{8}(4b_2-3)  +  \frac{3}{2}b_2b_3+ \frac{1}{8}(13b_2 -12) - \frac{21}{8}= b_6 - b_5 + \frac{3}{2}b_4 - \frac{5}{2}b_3 + 24b_2 + b_2b_3 -\frac{1}{2}b_2b_4-\frac{1}{2}= $.

$b_6 -b_5 + \frac{3}{2}b_4 - \frac{5}{2}b_3 - \frac{1}{2}b_2 + b_2b_3 - \frac{1}{2}b_2b_4 + \frac{1}{2}$.

\

\noi The computations can be done {\it \ˆ la main}; they are faster on a computer. We recover the formulas stated earlier, in the introduction, in particular,
$$2d_5(b _2,b _3,b _4,b _5) = 2b _5 - 3 b _4  + 5 b _3 + b _2+ b _4b _2 - 2 b _3b _2  - 1,$$
$$4d_6 = 4b _6 - 6 b _5+ 10 b _4 - 13 b _3+ 6 b _2 + 2b _5b _2  + 2 b _4b _3  - 6 b _4 b _2 + 3 b _3 b _2   - 6.$$

\

\noi One step further, we get expressions of those formulas with the $c_n$'s instead of the  $b_n$'s, using the relations\nts:
\[2c_{2n+1} = b_{2n+1} \anda 2c_{2n} = b_{2n} + a_n.\tag{$**$}\]
Thus,

$2a_2  = b_2 - 1$ 

$2a_3  = b_3 -  b_2 + 1$ 

$2a_4 =  b_4 - b_3+\frac{1}{2}b_2 -\frac{1}{2}$
  
$2a_5 = b_5 -b_4 -\frac{1}{2}b_2b_3+\frac{3}{2}b_3 +  \frac{1}{2}b_2- \frac{1}{2}$ 

\noi become, respectively,

$a_2 = c_2 -1$,

$a_3 = c_3 - c_2 + 1$,

$a_4 = c_4 - c_3$,
 
$a_5  = c_5 - c_4 + 2c_3 + c_2 - c_2c_3 - 1$.

\

{\bf As it clearly appears, now,  the  expressions of the $\bf a_n$'s as functions of the $\bf c_k$'s are  simpler than the preceding expressions as functions  of the $\bf b_k$'s.} Indeed,  those are polynomials in the $c_k$'s  all of whose coefficients are integers.

\

\noi Whence the question\nts: Can $a_n$ always be presented as a  polynomial from  the ring $\Z[c_1,c_2,\dots]$\nts? The answer is yes\! ! As we shall presently show, after having established the following more general result.

\

\

\

\

\head{Reciprocal forms}

\

For brevity's sake,  we use the following neologism (till the end of the paper, only)\nts: A {\bf compliform polynomial} is a polynomial whose  coefficients are integers and has degree $1$ in each of the variables.

\

\noi Similarly, in order to be short, and all along the text, we shall denote  M the algebra of polynomials in $x_2,x_3,\dots$, having integers for coefficients, and denote C the set of   compliform polynomials in $x_2,x_3,\dots$, so that we have C $\inc$ M $= \Z[x_2,x_3,\dots]$.

\su{Proposition} \sl In any given commutative ring $R$ (with identity) take two infinite sequences, $y_2,y_3,\dots,$and $z_2,z_3,\dots,$ where the $y_n$'s are  idempotent elements,  and such that\nts :
$$z_{m+1} = y_{m+1} + s_m(y_2,\dots,y_m)$$
where the $s_m$'s are  polynomials having integers for coefficients\rm.

\

\noi (i) \sl We, then, have\nts :
$$y_{n+1} = z_{n+1}  + t_n(z_2,\dots,z_n)$$
where the $t_n$'s are   compliform polynomials\rm.

\

\noi (ii) \sl Moreover, for each polynomial, $u(x_2,\dots,x_n)$, having integers for coefficients,  there is a compliform polynomial, $v(x_2,\dots,x_n)$,  such that\nts :
$$u(y_2,\dots,y_n)= v(z_2,\dots,z_n).$$\rm

\su{Proof} Proceed by induction. Suppose that, for a given $m\geqs 2$ and for each index $k$,  $2\leqs k \leqs m$,   the  element $y_k$ is written  as a compliform polynomial $y_k = z_k + t_{k-1}(z_2,\dots,z_{k-1})$. This is true for $m=2$ since $y_2 = z_2 - s_1$ and $s_1$ is an integer. Write  $y_k^2 = y_k$ and get  an expression of  $z_k^2$
  as  a polynomial in $z_2,\dots,z_k$, having integers for coefficients, and degree $1$ in $z_k$.  Using the algorithm above,   each polynomial in $z_2,\dots,z_k$, with  integers for coefficients,  can be  {\bf reduced}to a compliform polynomial.
Proceed and reduce the polynomial 
$$ s_m(\dots, z_k +  t_{k-1}(z_2,\dots,z_{k-1}), \dots)=  s_m(y_2,\dots,y_m)= z_{m+1} - y_{m+1}.$$
The expression thus obtained
ends the induction, and the proof of the first part (i). To prove (ii), write
$$u(y_2,\dots,y_n) = u(z_2 + t_1, \dots, z_n + t_{n-1}(z_2,\dots,z_n))$$
then reduce this last polynomial  to compliform.
\qed

\

\su{Corollary}{\sl For each integer $n\geqs 1$, there is a compliform polynomial, $e_n(x_2,\dots,x_n)$, such that,  in the algebra ETB,  we have\nts :
\[c_{n+1} = a_{n+1} + e_n(c_2,\dots,c_n).\tag{7}\]}
Moreover, there is a simple algorithm to compute those polynomials.

\

\noi Indeed, just apply the proposition, with the identities
\[c_n = \sum_{0\leqs k\leqs n/2} a_ka_{n-k},\tag{5}\]
written as
$$c_{m+1} = a_{m+1} + s_m(a_2,\dots,a_m).$$

\

\noi The proposition does not apply to   $b_n$ because
$b_n = 2a_n + \sum_{1\leqs k < n} a_ka_{n-k}$ is not of the required form.

\

\su{A first very special case} Getting on with the previous computations of  $a_2,a_3,a_4,a_5$, we reach the following\nts: 
$$a_6 = c_6 - c_5 + 2c_4 - 4c_3 + 2c_3c_2 - c_4c_2.$$
 A small surprise is  concealed herein\nts: Setting $c_2 = c_3 = c_4 = c_5 = c_6 = 1$, one gets $a_6 = -1$ which is not  compatible with a basis\! ! We state this result as follows\nts: Let  $p_2, p_3, p_4, p_5,p_6,$ be the  numbers  of $A$-presentations associated to any given basis $A$; then  at least one of those numbers is $\geqs 2$\! ! Of course, it follows that one, at least, of $r_2,r_3,r_4,r_5,r_6,$ is $\geqs 3$, due to the relations between the $r_n$'s and the $p_n$'s.

\

\noi In the same vein,  in 2003, {\sc Grekos, Haddad, Helou}, and {\sc Pihko} [3] showed that, for each basis, one at least of the representation numbers $r_n$ is $\geqs 6$.    The result was further improved  in 2006, by {\sc Borwein, Choi}, and {\sc Chu} [1], to $r_n\geqs 8$,

\

\head{Compliform polynomials\\ A parenthesis}

\

Look again at M $= \Z[x_2,x_3,\dots]$, the algebra of all polynomials in the variables $x_2,x_3,\dots$,  and integers for coeffcicients. It is the free (communtative) algebra generated by the $x_n$'s over $\Z$.  We have  denoted C the set of all compliform polynomials with the $x_n$'s for variables.  This is a subset of the free algebra M. It is, by no means, a subalgebra. However, it is a submodule, a subgroup with respect to addition.

\

\noi Let us elaborate. To each  polynomial $u(x_2,x_3,\dots,x_n)\in$ M is attached the compliform polynomial $\hat u(x_2,x_3,\dots,x_n)$ obtained on  substituting in $u$ an $x_k$ to each monomial $x_k^j$, for $j\geqs 2$. The map $u\mapsto \hat u$ from M to C  is an {\it onto} group homomorphism\nts : It is a {\bf retraction} of the additive group M onto its subgroup C.

\su{Scholium} \sl A  compliform polynomial, $u(x,y,\dots,z)$, on $n$ variables is identically zero if, for some  choice of given couples of values  of the variables, $x_0 \neq x_1, y_0\neq y_1, \dots, z_0 \neq z_1$, all the values $u(x_i,y_j,\dots,z_k)$ of the polynomial are zero\rm.

\

\noi {\bf Indeed}, it is clear that a compliform polynomial on one variable is identically zero if its values are zero for  two different values of the variable. Also,  a  compliform polynomial  $u(x,y)$ on two variables is identically zero if, for given values $x_0 \neq x_1$ and $y_0\neq y_1$, we have $u(x_i,y_j) = 0$. We just have to write $u(x,y) = v(y)x + w(y)$
with $u(y)$ and $w(y)$, both, compliform. Then notice that the two polynomials  on one variable, $u(x,y_0)$ and $u(x,y_1)$, are zero, so that  $v(y_0), v(y_1), w(y_0), w(y_1)$ are zero and, therefore, $v(y)$ and $w(y)$ are identically zero. We  conclude  by induction. {\bf qed}

\

\noi \sl For any $u(x_2,x_3,\dots,x_n)\in$ M, the following  statements are equivalent\rm\nts:  

\

(i) \ \  \sl $u(x_2,x_3,\dots,x_n)= 0$ whatever $x_k\in \{0,1\}$.\rm

(ii) \ \sl $\hat u(x_2,x_3,\dots,x_n)\equiv 0$.\rm

\

\noi Indeed, if $x_k\in \{0,1\}$, then $u(x_2,x_3,\dots,x_n)=  \hat u(x_2,x_3,\dots,x_n)$. What is left follows, immediately, from the preceding result.

\

\noi We also notice that a compliform  polynomial  $u(x,y,\dots,z)$ is an affine function of each of its variables. For instance, we have\nts:
$$u(s+t,y,\dots,z) = u(s,y,\dots,z)+ u(t,y,\dots,z)- u(0,y,\dots,z).$$

\

\head{One last avatar\\
The reduced Erd\šs-Tur\'an algebra}

\

The kind of reciprocity  highlighted above, as well as the lemma, in particular, show how to dispense, reasonably, with fractional coefficients in the algebra ETB. This leads to the following ultimate simplification.

\

\noi We call {\bf reduced  Erd\šs-Tur\'an algebra}, and denote ETR, the reduced subalgebra of  ETB generated  by $\Z$ and the $a_n$'s. That simply is $\Z[a_2,a_3,\dots]$, the algebra of polynomials,  endowed with the idempotency relations\nts
: \[a_n^2 = a_n.\tag{$***$}\]
Notice that the $b_n$'s and $c_n$'s belong to ETR since we have\nts:
\[b_n = a_0a_n + a_1a_{n-1} + \dots + a_ka_{n-k} + \dots +
a_na_0,\tag{$*$}\]
\[c_n =  \sum_{0\leqs k \leqs n/2} a_ka_{n-k}.\tag{5}\]
We could have introduced ETR earlier, directly, without all the detours. We rather chose to follow the path that led to it\! !
The algebra ETR has quite a few different, equivalent,  presentations. For instance, ETR is the algebra of polynomials  $\Z[a_2,a_3,\dots;b_2,b_3,\dots;c_2,c_3,\dots]$ with the relations ($*$), ($**$), and ($***$). It also is the algebra of polynomials $\Z[a_2,a_3,\dots;c_2,c_3,\dots]$ with the relations  ($***$) and (5) in which the $b_n$'s would be defined by  ($*$). 

\

\noi Still, the reduced Erd\šs-Tur\'an  algebra is also the quotient of the free algebra M$=\Z[x_2,x_3,\dots]$  by its ideal $I$ generated by the $(x_n^2-x_n)$'s. So to say, it is the algebra of polynomials, with integers for coefficients, in the {\bf idempotent} variables $a_2,a_3,\dots$. Notice that each non zero element in the ideal  $I$ is a polynomial with, at least, one variable having degree $\geqs 2$.

\

\head{On the structure of the algebra ETR}

\

Of course, all the relations ($*$), ($**$), ($***$), (5), and (7),  still hold in the algebra ETR.

\

\su{S1 How the $\mathbf {a_k}$'s generate ETR} The kernel of the following map from C to ETR\nts :
$$u(x_2,x_3,\dots,x_n) \mapsto u(a_2,a_3,\dots,a_n)$$
 is zero since $u$ has degree 1 at most in each of the variables. This map, therefore, is  a group isomorphism, and we have\nts:
$$\mathrm{ETR} = \{u(a_2,a_3,\dots,a_n) : u\in \mathrm C\}.$$
Thus, for $u(x_2,\dots,x_n)\in$ C, we have $u(a_2,a_3,\dots,a_n) = 0$ if and only if  $u(x_2,\dots,x_n) = 0$.

\

\noi For $u(x_2,x_3\dots x_n) \in$ M, we clearly have  $u(a_2,a_3,\dots,a_n) = \hat u(a_2,a_3,\dots,a_n)$.

\

\noi It should be remembered that, to each element $z$ of the algebra ETR is associated a unique compliform polynomial  $q_z = q_z(x_2,x_3,\dots,x_n)\in$ C, the one for which $z = q_z(a_2,a_3,\dots)$. Also associated to $z$ is the set of all polynomials  $u(x_2,x_3,\dots,x_n)\in$ M such that $\hat u = q_z$. 

\su{S2 Universality of ETR} Given a commutative ring (with identity) $R$ and an infinite sequence of idempotent elements, $i_2, i_3, \dots$, in $R$, the map $a_n\mapsto i_n$, for $n\geqs 2$, has a unique extension as a homomorphism of the ring  ETR into the ring $R$\nts : The extension is the map  $z= q_z(a_2,\dots,a_n) \mapsto q_z(i_2,\dots,i_n)$. So, in a way, the algebra ETR is {\it universal} in the category of commutative rings generated by a (finite or infinite) sequence of idempotent elements.

\su{S3 Non integrity} The algebra ETR is far from being an integral domain. It has many divisors of zero. For instance, we have  $a_n(a_n-1) = 0$ though $a_n \neq 0$ and $a_n\neq 1$, for $n\geqs 2$. We also have $(c_2 -1)(c_2-2) = 0$ but $c_2\neq 1$ and $c_2\neq 2$. Quite generally, we have  $(c_{n+1} - e_n)(c_{n+1} - e_n-1) = 0$ while none of the two factors is zero.

\

\noi Here again is an  identity, less visible. We have $c_3 = a_3 + a_2$, therefore, $c_3(c_3-1)(c_3-2) = 0$, all computations done.

\

\head{Topology for an interlude}

\

Among all  bases of $\N$, we differentiate between those for which the sequence of  $r_n$'s is bounded  and the others. The first are the bases for which the Erd\šs-Tur\'an conjecture  does not hold. It can be shown, with not much effort, that if any of the first exist, their set is not very big, it is {\bf meager} (a first category set) in the  sense of Baire's theory. \{As is known, in topology, the notion of meager subsets is the counterpart of  negligible subsets  in measure theory.\}

\

\noi Let $S$ be the set of all  formal series having the form
$$u(t) = \sum u_nt^n \ \ \text{where $u_n \in \{0,1\}$ for eacht $n$}.$$ 
The set  $S$ is identified to the compact space  $K = 2^\N$, when identiying the sequence $(u_n)_{n\in\N}$ to the corresponding point  in $K$.  To each series $u \in K$, we associate its square $v(t) = u(t)^2 = \sum r_n t^n$. Let $T$ be the subspace in $K$ of all the series $u$ such that $r_n \geqs 1$ for each index $n$, and let $U$ be the  subspace in $T$ of all the $u \in T$ such that the sequence of $r_n$'s is bounded. Thus, quite specifically, one can think of  $T$ as the space of bases for $\N$ while $U$ is the subspace of those bases with a bounded profile (if any). We have a straightfordward proof that $U$ is a  \sl meager \rm subset in $T$. 

\

\noi This means that {\bf $\mathbf U$ is contained in a countable union of closed subsets with empty interiors in $\mathbf T$.}

\su{Proof} Set $V = T \setminus U$. For each integer $p>0$, set
 $$V(p) = \{u \in T : \exists m  \ (r_m \geqs p)\}.$$
Thus, $V = \bigcap V(p)$. On the one hand, each $V(p)$ is open in $T$\nts: Indeed, given $u \in V(p)$, take an index $m$ such that $r_m \geqs p$; each series $h \in T$ whose coefficients coincide with those of $u$ up to rank $m$ also belongs to $V(p)$. On the other hand, $V$ is dense in $T$\nts: Indeed, given $h\in T$ and any index $q$, let  $u$ be a series whose first $q$ coefficients coincide with those of $h$, and all other coefficients are equal to 1\nts: The series $u$ belongs to $V$. So, $V$ is a \sl residual \rm (second category) in $T$ and its complement is, therefore,  \sl meager\rm.\qed

\

\noi{\bf In the set of all bases, the ones that satisfy the Erd\šs-Tur\'an conjecture are \emph{generic}}, i.e., their set is co-meager.

\

\noi \{That should explain why, if the  Erd\šs-Tur\'an conjecture is false, it is difficult  to find a counterexample. This is reminiscent of the analogous  situation  when people were looking for real continuous functions having nowhere, neither finite nor infinite, derivatives, in the   20's  of the past century; read in {\sc Saks} [10] about that.\}

\

\head{The main system}

\

\noi Keeping the same notations,  pretend the relations (7) of the lemma are a system of equations with the $a_n$'s as given data, and the $c_n$'s as unknowns, [not much confusion to fear]. So, write\nts:
\[c_{n+1} = a_{n+1} + e_n(c_2,\dots,c_n).\tag{SP}\]
Start with the first equations, and get

\

$c_2 = a_2 + 1$,

$c_3 = a_3 + c_2 - 1$,

$c_4 = a_4 + c_3$,

$c_5 = a_5 + c_4 - 2c_3 - c_2+ c_3c_2 +1$,

$c_6 = a_6 +c_5 - 2c_4+ 4c_3+ c_4c_2 - 2c_3c_2$,

$c_7 = a_7  + c_6 - 2 c_5  + 4 c_4 - 4 c_3 + 3 c_2 +  c_5c_2 + c_4 c_3 - 3 c_4 c_2  + c_3 c_2  
 - 3$,

$c_8= a_8 + c_7- 2 c_6+ 4 c_5 - 4 c_4 + 4 c_3  - 4 c_2+ c_6 c_2 + c_5 c_3- 3  c_5c_2- 2 c_4 c_3$,
 
\qquad$+4 c_4 c_2 - c_3 c_2 +4 $,

$\dots\dots\dots\dots\dots\dots\dots$

$c_{n+1} = a_{n+1} + e_n(c_2,\dots,c_n)$,

$\dots\dots\dots\dots\dots\dots\dots$

\

\noi Looking  at its structure, taking  into account, solely, the fact that the polynomials $e_n(c_2,\dots,c_n)$ are compliform, many properties of the system SP can be dervived. Here are the first two of them.
\su{P1 Existence of the solutions} Whatever the values given to the sequence  $a=(a_2, a_3, \dots, a_n\dots)$ are, the system has a unique solution  $c[a] =(c_2, c_3, \dots, c_n\dots)$. If the values in $a$ are real numbers  (resp., complex numbers), then so are the values in $c[a]$. If the values in  $a$ are integers, the values in $c[a]$ will also be integers. Quite generally, if the $a_k$'s belong to a given commutative ring (with unity),  the  $c_k$'s will also belong to that ring. 

\su{P2 Uniqueness of the system} Let
\[c_{n+1} = a_{n+1} + t_n(c_2,\dots,c_n),\]
be another system with compliform polynomials $t_n(c_2,\dots,c_n)$, again. If the two systems have the same solutions, then each polynomial $t_n$ is identical to the corresponding $e_n$. In fact, it is enough to  suppose  that, for each of the choices  $a_n\in\{0,1\}$, the two systems have the same solution $c[a]$, in the ring $\Z$, to reach the same conclusion\nts: The $t_n$'s are identical to the corresponding $e_n$'s. This is an easy consequence of the scholium, above.

\

Beyond its sole structure, the system SP has other noteworthy properties deriving from the way it was constructed.

\

\su{P3 Presentations}  Le $A$ be any subset of $\N$ such that $0$ and $1$ belong to it. Go back to the algebra ETR and set $a_n = 1 \ou 0$, according as $n$ belongs or not to $A$. Then, the solution in $c= (c_2,c_3,\dots)$ of the  system SP corresponding to the data $a= (a_2,a_3,\dots)$ is precisely $c_n = p_n$, the number of  $A$-presentations of $n$. That is how the system was constructed, anyway.

\su{Let us insist}At each step, according as $n+1$ belongs or not to $A$, we have\nts:
\[p_{n+1} = e_n(p_2,\dots,p_n) +1 \ora p_{n+1} = e_n(p_2,\dots,p_n).\tag{8}\]
For each index, there is a choice between two different values for $p_{n+1}$\nts: One of the two integers  $x_0 = e_n(p_2,\dots,p_n)\geqs 0$ or $x_1= x_0 + 1$. 
 One can imagine tossing it {\it heads or tails}. The infinite sequence of $p_n$'s is  thus generated {\it almost} automatically. Now, the word automatic having already been used in a precise meaning, to describe infinite sequences generated by special automata,  we shall use another word to describe the infinite sequences $1,p_2,p_3,\dots,$ associated to the subsets $A\inc \N$ (to which $0$ and $1$ belong)\nts : We call them ({\it faute de mieux}) the {\bf hemitropic} sequences. 

\

\noi We further observe the following. The subset $A$ is a basis of  $\N$ if and only if $p_n\geqs1$ holds, for each  $n$. To make it short, we call a hemitropic sequence $1,p_2,p_3,\dots,$  {\bf basic} when $p_n \geqs 1$ holds, for each $n\geqs1$.

\

\noi Erd\šs-Tur\'an's conjecture just says that {\bf no  hemitropic basic sequence is bounded}.
\su{A small digression} Hemitropic bounded sequences exist, some of them are even bounded by $1$, and not terminating by $0$'s. Those are associated to Sidon sets. Take $A=\{0,1,3,\dots,\}= \{v_{n+1} = 2v_n + 1 : v_0 = 0, n\in \N\}$  as a special case\nts: The  corresponding $p_n$'s are, all, either $0$ or $1$.

\

\head{Measure for measure}

\

Let us take a break. It can be shown that, among bases of $\N$, it is quite unlikely to find one that contradicts the Erd\šs-Tur\'an conjecture. In a private communication, 
{\sc Jean-Pierre KAHANE} showed me very precisely how to establish that type of result.

\

\noi Take any real number $p\in \ ]0,1]$. Set $P(a_n=1) = p$, $P(a_n = 0) = 1-p = q$. This means that, in a draw to determine $a_n\in \{0,1\}$, the probability to get a $1$ is equal to $p$. This also means that the discrete set  $\{0,1\}$ is endowed with a probability measure $P$ such that  $P(\{1\}) = p$ and $P(\{0\})= q$, then,  the product space  $K = \{0,1\}^\N = 2^\N$ is endowed with the product measure. Introduce the following random variables\nts:
$$c_n  = \sum_{0\leqs i\leqs n/2} a_ia_{n-i}= a_0a_n + a_1a_{n-1} + \cdots.$$
We show that, for any given integer $k$, the following holds\nts: 
$$P\left(\sup_{n\geqs 0} c_n \leqs k\right) =0.$$
This means that the sequence of $c_n$'s is  unbounded, {\bf almost surely}. So, the set of bases of $\N$ which contradict  the Erd\šs-Tur\'an conjecture is  a {\bf null set}, whatever the probability  $p>0$ is.

\su{Proof}  The event  $c_{2m+1} = k$ occurs if and only if there is a $k$ element subset $I\inc \{0,1,\dots,m\}$  such that
$$a_ia_{2m+1-i} = 1, \ \text{for $i\in I$, and $a_ia_{2m+1-i} = 0$, oherwise}.$$
This means that $I = \{i : a_i = a_{2m+1-i} = 1 \et i\leqs m\}$ is a $k$ element set. Therefore\nts:
$$P(c_{2m+1} = k) = \binom{m+1} k p^{2k}(1-p^2)^{m+1-k}.$$
Likewise,  the event  $c_{2m} = k$ occurs if and only if there is a $k$ element subset $I \inc \{0,1,\dots,m\}$ such that
$$a_ia_{2m-i} = 1, \ \text{for $i\in I$, and $a_ia_{2m-i} = 0$, otherwise}.$$
There are two possibilities, according as $m$ belongs to $I$ or not, therefore
$$P(c_{2m} = k) = \binom m {k-1} p^{2k-1}(1-p^2)^{m+1-k} + q\binom m k p^{2k}(1-p^2)^{m-k}.$$
By Stirling's formula, for any given $k$, the limit of $P(c_n=k)$ is $0$, when $n$ goes to infinity. The result follows.\qed

\

\noi In the same vein, going a bit further into computations, it can be shown that, whatever the probability $p > 1/2$ is, the following holds, almost surely\nts:
$$\limsup_{n} \frac{8c_n}{n}\geqs 1.$$

\

\head{A pause}

\

We have seen, above, why and how topology leaves very little space to  bases that would possibly contradict the Erd\šs-Tur\'an conjecture. Just now, we see that probability gives those bases no chance at all. Thus, the search for a counterexample to the conjecture {\bf \emph{ET}}  might  very much look like  {\sl walking on eggshells}, or a passage through {\sl many narrow doors}, a journey on {\sl a crest road} edged with precipices.

\

\noi Can we learn more from Algebra\! ? From Combinatorics\! ?

\

\head{The binary tree of hemitropic sequences} 

\

The set of all hemitropic sequences can be viewed as a tree $\cal T$. The nodes of the
tree $\cal T$ are the finite sequences $(1,p_2,\dots,p_n)$  where $(1,p_2,p_3,\dots)$ is a hemitropic sequence.
The root of the tree is the one term sequence $(p_1=1)$. Each
node $(1,p_2,\dots,p_n)$ has exactly
two immediate successors, a
{\bf lower node}, $(1,p_2,\dots,p_n,x)$, 
and an
{\bf upper node}, 
$(1,p_2,\dots,p_n,x+1)$. The
(infinite) branches of the tree are the hemitropic sequences $(1,p_2,\dots,p_n,\dots)$.

\

\noi A node (resp. branch) will be called \bf persistent \rm when
none of its $p_n$'s is $0$. In the tree $\cal T$, there is
always, at least, one  persitent branch passing through
any given persistent node\nts : Indeed, after each node, taking the immediate upper successor obtains a persistent branch.

\

\noi To each 
branch $p=(1,p_2, \dots)$ corresponds the subset $A \inc \N$ such that $p_n$ is the number of $A$-presentations of $n$. The branch $p$ is persistent if and only if $A$ is a
basis of $\N$\nts: This establishes a one-to-one  relation
between   persistent branches and bases.

\

\noi The conjecture {\bf \emph{ET}}can be stated as follows\nts:
The tree $\cal T$ has no persitent branch  $p=(1,p_2, \dots)$ where the set of $p_n$'s is bounded.

\

\

\

\head{As a tentative conclusion}

\

In a very precise meaning, we say that there is little room and almost no chance for the Erd\šs-Tur\'an conjecture to be wrong. Little room\nts? Or no room at all\nts!

\

\centerline{To be continued}

\

\noi [Maybe, just in case, I should apologize for \lq\lq bad English", here and there.]

\

\

\

\

\head{Remerciements\\ Acknowledgments} 

\

Tout d'abord, je dois remercier mes coll\gues et amis, Georges Grekos, Charles Helou et Jukka Pihko, avec lesquels j'ai de  nombreuses conversations instructives depuis plus d'une douzaine d'ann\Žes. Je tiens \ˆ remercier \Žgalement Monsieur le Professeur Jean-Pierre Kahane pour m'avoir appris, dans le d\Žtail, comment utiliser les probabilit\Žs dans ce contexte arithm\Žtique. Je remercie aussi toutes les  personnes qui ont \oe uvr\Ž \ˆ  la construction du logiciel Maple : il permet \ˆ tout un chacun d'\Žgaler les grands, comme Euler ou Gauss, en mani\re de calcul. Je remercie enfin tous ceux qui ont con\c cu et r\Žalis\Ž le site {\tt arXiv} lequel rend de signal\Žs services \ˆ la communaut\Ž des math\Žmaticiens.

\

First of all, I would like to thank my colleagues and friends, Georges Grekos, Charles Helou, and Jukka Pihko, for the  enlightening conversations I have with them, since more than twelve years, now. I also want to thank Jean-Pierre Kahane who showed me, in great detail, how to use probability  in this number theoretic context. I would like, also, to thank all those who worked on Maple, a software that makes it so easy to emulate the great, in  computation matters.  I also thank, very much indeed, those who created the site {\tt arXiv} providing a quite remarkable service to all mathematicians.

\

\head{Bibliographie\\ References}

\

\

[1] Borwein, Peter; Choi, Stephen; Chu, Frank.
An old conjecture of Erd\šs-Tur\'an on additive bases. Math. Comput. 75, no. 253, (2006) 475-484.

\

[2] Erd\šs, P.; Tur\'an, P. On a problem of Sidon in additive number theory, and on some related problems. J. London Math. Soc. 16 (1941)  212--215. 

\

[3] Grekos, G.; Haddad, L.; Helou, C.; Pihko, J.
On the Erd\šs-Tur\'an conjecture. 
J. Number Theory 102, no. 2, (2003)  339-352 .

\

[4] Grekos, G.; Haddad, L.; Helou, C.; Pihko, J.
Additive bases representations and the Erd\šs-Tur\'an conjecture. 
Chudnovsky, David (ed.) et al., Number theory: New York seminar 2003. New York, NY, (2004) 193-200 .

\

[5] Grekos, G.; Haddad, L.; Helou, C.; Pihko, J.
The class of Erd\šs-Tur\'an sets. 
Acta Arith. 117, no. 1, (2005)  81-105.

\

[6] Grekos, G.; Haddad, L.; Helou, C.; Pihko, J.
Variations on a theme of Cassels for additive bases. 
Int. J. Number Theory 2, no. 2,  (2006) 249-265.

\

[7] Grekos, G.; Haddad, L.; Helou, C.; Pihko, J.
Representation functions, Sidon sets and bases. 
Acta Arith. 130, no. 2, (2007) 149-156.

\

[8] Haddad, L.; Helou, C.; Pihko, J.
Analytic Erd\šs-Tur\'an conjectures and Erd\šs-Fuchs theorem. 
Int. J. Math. Math. Sci. 2005, no. 23, (2005) 3767-3780 .

\

[9] Nathanson, Melvyn B. Every function is the representation function of an additive basis for the integers. Port. Math. (N.S.) 62 (2005), no. 1, 55--72. 

\

[10] Saks, S. On the functions of Besicovitch in the space of continuous functions. Fundam. Math. 19 (1932) 211-219 .

\

[11] S\'andor, Csaba.
A note on a conjecture of Erd\šs-Tur\'an. 
Integers 8, no. 1, (2008) Article A30.

\

\

\

\enddocument